\documentclass[12pt]{article}

\usepackage{amssymb}
\usepackage{latexsym}
\usepackage{calc}
\usepackage{datetime}
\usepackage{graphicx}
\usepackage{amsmath,amsfonts}
\usepackage[framemethod=tikz]{mdframed}
\usepackage{color, colortbl}

\newcounter{hours}\newcounter{minutes}

\def\dom{{\rm dom \,}}
\def\beq{\begin{equation}}
\def\eeq{\end{equation}}

\numberwithin{equation}{section}

\newtheorem{theorem}{Theorem}
\numberwithin{theorem}{section} 
\newtheorem{lemma}{Lemma}
\numberwithin{lemma}{section} 
\newtheorem{corollary}{Corollary}
\numberwithin{corollary}{section} 

\newtheorem{assumption}{Assumption}
\newtheorem{definition}{Definition}

\newtheorem{example}{Example}
\newtheorem{remark}{Remark}
\numberwithin{remark}{section} 

\def\ba{\begin{array}}
\def\ea{\end{array}}
\def\beann{\begin{eqnarray*}}
\def\eeann{\end{eqnarray*}}
\def\bea{\begin{eqnarray}}
\def\eea{\end{eqnarray}}

\def\BT{\begin{theorem}}
\def\ET{\end{theorem}}
\def\BL{\begin{lemma}}
\def\EL{\end{lemma}}
\def\BC{\begin{corollary}}
\def\EC{\end{corollary}}
\def\BE{\begin{example}}
\def\EE{\end{example}}
\def\BD{\begin{definition}}
\def\ED{\end{definition}}
\def\BR{\begin{remark}}
\def\ER{\end{remark}}
\def\BAS{\begin{assumption}}
\def\EAS{\end{assumption}}
\def\BI{\begin{itemize}}
\def\EI{\end{itemize}}

\def\BMP{\begin{minipage}{9.5cm}}
\def\EMP{\end{minipage}}
\def\MPT{\begin{minipage}{11.5cm}}
\def\EPT{\end{minipage}}

\def\la{\langle}
\def\ra{\rangle}

\title{\vspace{10mm} \textbf{Adaptive Third-Order Methods for Composite Convex Optimization}}

\author{Geovani Nunes Grapiglia\thanks{Université catholique de Louvain, ICTEAM/INMA, Avenue Gerges Lema\^{i}tre 4-6, 1348 Louvain-la-Neuve, Belgium (geovani.grapiglia@uclouvain.be). This author was partially supported by the National Council for Scientific and Technological Development - Brazil (grant 312777/2020-5) and by the European Research Council Advanced Grant 788368.} 
        \and Yurii Nesterov\thanks{Université catholique de Louvain, Center for Operations Research and Econometrics (CORE), 34 voie du Roman Pays, 1348 Louvain-la-Neuve, Belgium (Yurii.Nesterov@uclouvain.be). This author was financed by the European Research Council Advanced Grant 788368.}}

\begin{document}

\maketitle

\maketitle

\abstract{In this paper we propose third-order methods for composite convex optimization problems in which the smooth part is a three-times continuously differentiable function with Lipschitz continuous third-order derivatives. The methods are adaptive in the sense that they do not require the knowledge of the Lipschitz constant. Trial points are computed by the inexact minimization of models that consist in the nonsmooth part of the objective plus a quartic regularization of third-order Taylor polynomial of the smooth part. Specifically, approximate solutions of the auxiliary problems are obtained by using a Bregman gradient method as inner solver. Different from existing adaptive approaches for high-order methods, in our new schemes the regularization parameters are adjusted taking into account the progress of the inner solver. With this technique, we show that the basic method finds an $\epsilon$-approximate minimizer of the objective function performing at most $\mathcal{O}\left(|\log(\epsilon)|\epsilon^{-\frac{1}{3}}\right)$ iterations of the inner solver. An accelerated adaptive third-order method is also presented with total inner iteration complexity of $\mathcal{O}\left(|\log(\epsilon)|\epsilon^{-\frac{1}{4}}\right)$.}

\vspace{10ex}\noindent


\thispagestyle{empty}

\newpage\setcounter{page}{1}

\section{Introduction}
\setcounter{equation}{0}

\subsection{Motivation}

Following \cite{NP}, several methods based on regularized Taylor models have been proposed for unconstrained minimization of convex and nonconvex functions with Lipschitz continuous $p$th derivatives (see, e.g., \cite{BAES,Birgin,CGT2,DOI,DU,GAS,GN0,MR}). In these methods, trial points are obtained by the minimization of $(p+1)$-regularizations of the $p$th-order Taylor approximations of the objective function ($p\geq 2$). To ensure fast global convergence rates, it is enough to compute suitable approximate stationary points of the models \cite{Birgin}. Moreover, if the objective function is convex, then any $(p+1)$-regularization of the corresponding $p$th-order Taylor approximation is also convex when the regularization parameter is sufficiently large \cite{NES1}. Remarkably, for $p=3$, additional relative smoothness properties \cite{BBT,LU} were proved in \cite{NES1} for third-order tensor models. These properties allow the use of Bregman gradient methods to approximately minimize third-order tensor models with linear rate of convergence \cite{NES1}. 

The fast minimization of third-order tensor models can be guaranteed only when the regularization parameter is bigger than some multiple of the Lipschitz constant of the third-order derivative of the objective function. The simplest way to satisfy this condition is to use explicitly such Lipschitz constant, when it is known. As pointed in \cite{NES1}, this requirement is one of the main limitations for practical implementations of third-order tensor methods, since in general Lipschitz constants of third-order derivatives are unknown. In this paper we address this problem. Specifically, we propose adaptive third-order methods for composite convex optimization problems. The me\-thods are adaptive in the sense that they do not require the knowledge of the Lipschitz constant. Trial points are computed using a Bregman gradient method as inner solver. We propose a novel way to update the regularization parameters in the third-order models used to compute trial points. Different from existing adaptive tensor methods \cite{GN01,JIANG}, in our new methods the regularization parameters are adjusted taking into account the progress of the inner solver. When the regularization parameter is sufficiently large, the inner solver is guaranteed to have linear rate of convergence. We design an implementable way to check the violation of linear convergence. Thus, once a slow convergence rate is detected, the execution of the inner solver is stopped and the regularization parameter is increased. With this technique, we show that our basic adaptive third-order method finds an $\epsilon$-approximate minimizer of the objective function performing at most $\mathcal{O}\left(|\log(\epsilon)|\epsilon^{-\frac{1}{3}}\right)$ iterations of the inner solver, while our accelerated adaptive third-order method has total inner iteration complexity of $\mathcal{O}\left(|\log(\epsilon)|\epsilon^{-\frac{1}{4}}\right)$. Preliminary numerical results are also reported and illustrate our new adaptive technique.

\subsection{Contents}
The paper is organized as follows. In Section 2 we define our problem and present the main results about the Bregman gradient method applied to the minimization of third-order tensor models. In Section 3 we present and analyse the basic adaptive third-order method. Section 4 is dedicated to the accelerated adaptive third-order method. Illustrative numerical results are reported in Section 5. Finally, in Section 6 we discuss the possibility of second-order implementations of our adaptive schemes.

\subsection{Notations}

In what follows, $\|\,\cdot\,\|$ denotes the Euclidian norm in $\mathbb{R}^{n}$, or the matrix and tensor norms induced by it. For a set $S\subset\mathbb{R}^{n}$, $\text{conv}(S)$ denotes its convex hull. Given a matrix $A\in\mathbb{R}^{n\times n}$, $\text{trace}(A)=\sum_{i=1}^{n}A_{ii}$. Moreover, for a smooth function $f:\mathbb{R}^{n}\to\mathbb{R}$, we denote by $D^{3}f(x)$ the third-order derivative of $f(\,\cdot\,)$ at point $x$, which is the tensor defined by 
\begin{equation*}
\left[D^{3}f(x)\right]_{ijk}=\dfrac{\partial^{3}f}{\partial x_{i}\partial x_{j}\partial x_{k}}(x),\quad i,j,k\in\left\{1,\ldots,n\right\}.
\end{equation*}

\section{Problem Statement and Auxiliary Results}

In this paper we consider methods for solving the following composite minimization problem
\begin{equation}
\min_{x\in\mathbb{R}^{n}}\,\tilde{f}(x)\equiv f(x)+\psi(x),
\label{eq:2.1}
\end{equation}
where $f:\mathbb{R}^{n}\to\mathbb{R}$ is convex and three times differentiable, while $\psi:\mathbb{R}^{n}\to\mathbb{R}\cup\left\{+\infty\right\}$ is a simple convex function. More specifically, we will assume that the third-order derivative of $f(\,\cdot\,)$ is $L_{f}$-Lipschitz continuous, i.e., 
\begin{equation}
\|D^{3}f(x)-D^{3}f(y)\|\leq L_{f}\|x-y\|,\quad\forall x,y\in\mathbb{R}^{n}.
\label{eq:2.2}
\end{equation}
Given $x\in\mathbb{R}^{n}$, the third-order Taylor polynomial of $f(\,\cdot\,)$ around $x$ is given by
\begin{equation}
\Phi_{x}(y)=f(x)+\sum_{i=1}^{3}\dfrac{1}{i!}D^{i}f(x)[y-x]^{i},\quad y\in\mathbb{R}^{n}.
\label{eq:2.3}
\end{equation}
It follows from (\ref{eq:2.2}) and (\ref{eq:2.3}) that
\begin{equation}
|f(y)-\Phi_{x}(y)|\leq\dfrac{L_{f}}{4!}\|y-x\|^{4},
\label{eq:2.4}
\end{equation}
\begin{equation}
\|\nabla f(y)-\nabla\Phi_{x}(y)\|\leq\dfrac{L_{f}}{3!}\|y-x\|^{3},
\label{eq:2.5}
\end{equation}
and
\begin{equation}
\|\nabla^{2}f(y)-\nabla^{2}\Phi_{x}(y)\|\leq\dfrac{L_{f}}{2!}\|y-x\|^{2}
\label{eq:2.6}
\end{equation}
Define $d_{4}:\mathbb{R}^{n}\to\mathbb{R}$ by
\begin{equation}
d_{4}(x)=\dfrac{1}{4}\|x\|^{4},
\label{eq:2.7}
\end{equation}
and consider the following model of $\tilde{f}(\,\cdot\,)$ around of $x$:
\begin{equation}
\tilde{\Omega}_{x,M}(y)=\Omega_{x,M}(y)+\psi(y),\quad y\in\dom\psi,
\label{eq:2.8}
\end{equation}
where 
\begin{equation}
\Omega_{x,M}(y)=\Phi_{x}(y)+\dfrac{M}{2}d_{4}(y-x),
\label{eq:2.9}
\end{equation}
with $M>0$ being an estimate of $L_{f}$. If $M\geq L_{f}$, then it follows from (\ref{eq:2.4}) that
\begin{eqnarray}
\tilde{f}(y)&\leq & \Phi_{x}(y)+\dfrac{4L_{f}}{4!}d_{4}(y-x)+\psi(y)\nonumber\\
                 & =  & \Phi_{x}(y)+\dfrac{M}{3!}d_{4}(y-x)+\psi(y)\nonumber\\
                 &\leq & \Phi_{x}(y)+\dfrac{M}{2!}d_{4}(y-x)+\psi(y)\nonumber\\
                 &  =  & \tilde{\Omega}_{x,M}(y).
\label{eq:2.10}
\end{eqnarray}
This remark motivates the following third-order method to solve (\ref{eq:2.1}):
\\[0.2cm]
\begin{mdframed}
\begin{equation}
x_{t+1}\in\arg\min_{y\in\mathbb{R}^{n}}\tilde{\Omega}_{x_{t},M}(y),\quad t\geq 0.
\label{eq:2.11}
\end{equation}
\end{mdframed}
\vspace{0.2cm}
Given, $x\in\mathbb{R}^{n}$, notice that $\Omega_{x,M}(\,\cdot\,)$ in (\ref{eq:2.9}) is a multivariate polynomial. Therefore, in general, $\tilde{\Omega}_{x,M}(\,\cdot\,)$ in (\ref{eq:2.8}) is a nonconvex function, which makes its exact minimization in (\ref{eq:2.11}) out of reach. Remarkably, if $M\geq L_{f}$, it follows from Theorem 1 in \cite{NES1} that $\tilde{\Omega}_{x,M}(\,\cdot\,)$ is convex, allowing the use of efficient methods from Convex Optimization to approximately solve the problem
\begin{equation}
\min_{y\in\mathbb{R}^{n}}\tilde{\Omega}_{x,M}(y).
\label{eq:2.12}
\end{equation}
In the present paper, we are interested in the computation of an inexact solution $x^{+}$ of (\ref{eq:2.12}) such that
\begin{equation}
\|\nabla \Omega_{x,M}(x^{+})+g_{\psi}(x^{+})\|\leq \dfrac{M}{6}\|x^{+}-x\|^{3}
\label{eq:2.13}
\end{equation}
for some $g_{\psi}(x^{+})\in\partial\psi(x^{+})$. Such point $x^{+}$ enjoys the following property which motivates the use of this type of inexact solution of (\ref{eq:2.12}) in the iterative minimization of $\tilde{f}(\,\cdot\,)$.
\begin{lemma}
\label{lem:2.1}
Suppose that $f:\mathbb{R}^{n}\to\mathbb{R}$ satisfies (\ref{eq:2.2}) and let $x^{+}$ be a point for which (\ref{eq:2.13}) holds. If $M\geq L_{f}$ then
\begin{equation}
\langle\nabla\tilde{f}(x^{+}),x-x^{+}\rangle\geq \dfrac{1}{6M^{\frac{1}{3}}}\|\nabla\tilde{f}(x^{+})\|^{\frac{4}{3}},
\label{eq:2.14}
\end{equation}
where $\nabla\tilde{f}(x^{+})=\nabla f(x^{+})+g_{\psi}(x^{+})$.
\end{lemma}

\begin{proof}
By (\ref{eq:2.5}), (\ref{eq:2.9}) and (\ref{eq:2.13}) we have
\begin{eqnarray*}
\|\nabla\tilde{f}(x^{+})+\dfrac{M}{2}\|x^{+}-x\|^{2}(x^{+}-x)\|&=&\|\nabla f(x^{+})-\nabla\Phi_{x}(x^{+})+\nabla\Omega_{x,M}(x^{+})+g_{\psi}(x^{+})\|\\
&\leq & \|\nabla f(x^{+})-\nabla\Phi_{x}(x^{+})\|+\|\nabla\Omega_{x,M}(x^{+})+g_{\psi}(x^{+})\|\\
&\leq &\dfrac{L_{f}+M}{6}\|x^{+}-x\|^{3}.
\end{eqnarray*}
Since $M\geq L_{f}$, we get
\begin{equation*}
\|\nabla\tilde{f}(x^{+})+\dfrac{M}{2}\|x^{+}-x\|^{2}(x^{+}-x)\|\leq \dfrac{M}{3}\|x^{+}-x\|^{3}.
\end{equation*}
Consequently,
\begin{eqnarray*}
\dfrac{M^{2}}{9}\|x^{+}-x\|^{6}&\geq &\|\nabla\tilde{f}(x^{+})+\dfrac{M}{2}\|x^{+}-x\|^{2}(x^{+}-x)\|^{2}\\
                                                        & =     &\|\nabla\tilde{f}(x^{+})\|^{2}+M\|x^{+}-x\|^{2}\langle\nabla\tilde{f}(x^{+}),x^{+}-x\rangle+\dfrac{M^{2}}{4}\|x^{+}-x\|^{6},
\end{eqnarray*}
which implies that
\begin{eqnarray}
\langle\nabla\tilde{f}(x^{+}),x-x^{+}\rangle&\geq&\dfrac{\|\nabla\tilde{f}(x^{+})\|^{2}}{M\|x^{+}-x\|^{2}}+M\left(\frac{1}{4}-\frac{1}{9}\right)\|x^{+}-x\|^{4}\nonumber\\
&\geq &\dfrac{5M}{36}\|x^{+}-x\|^{4}.
\label{eq:2.15}
\end{eqnarray}
On the other hand, by (\ref{eq:2.5}), (\ref{eq:2.9}) and (\ref{eq:2.13}) we also have
\begin{eqnarray*}
\|\nabla \tilde{f}(x^{+})\|&\leq & \|\nabla f(x^{+})-\Phi_{x}(x^{+})\|+\|\nabla\Phi_{x}(x^{+})-\nabla\Omega_{x,M}(x^{+})\|+\|\nabla\Omega_{x,M}(x^{+})+g_{\psi}(x^{+})\|\\
&\leq &\dfrac{L_{f}}{3!}\|x^{+}-x\|^{3}+\dfrac{M}{2}\|x^{+}-x\|^{3}+\dfrac{M}{6}\|x^{+}-x\|^{3}\\
& = & \dfrac{(L_{f}+3M+M)}{6}\|x^{+}-x\|^{3}\\
&\leq & \dfrac{L_{f}+4M}{6}\|x^{+}-x\|^{3}.
\end{eqnarray*}
Since $M\geq L_{f}$, we get
\begin{equation*}
\|\nabla\tilde{f}(x^{+})\|\leq\dfrac{5M}{6}\|x^{+}-x\|^{3},
\end{equation*}
which implies that
\begin{equation}
\|x^{+}-x\|^{3}\geq\dfrac{6}{5M}\|\nabla \tilde{f}(x^{+})\|.
\label{eq:2.16}
\end{equation}
Finally, combining (\ref{eq:2.15}) and (\ref{eq:2.16}) we obtain (\ref{eq:2.14}). 
\end{proof}

\noindent In what follows we have a sufficient condition for (\ref{eq:2.13}).

\begin{lemma}
\label{lem:2.2}
Suppose that $f:\mathbb{R}^{n}\to\mathbb{R}$ satisfies (\ref{eq:2.2}), and let $x\in\mathbb{R}^{n}$ and $M,\epsilon>0$. If a point $x^{+}$ satisfies
\begin{equation}
\|\nabla f(x^{+})+g_{\psi}(x^{+})\|\geq\epsilon
\label{eq:2.17}
\end{equation}
and
\begin{equation}
\|\nabla\Omega_{x,M}(x^{+})+g_{\psi}(x^{+})\|\leq\dfrac{1}{2}\min\left\{1,\dfrac{M}{L_{f}+3M}\right\}\epsilon,
\label{eq:2.18}
\end{equation}
then $x^{+}$ satisfies (\ref{eq:2.13}).
\end{lemma}

\begin{proof}
By (\ref{eq:2.17}), (\ref{eq:2.5}) and (\ref{eq:2.18}) we get
\begin{eqnarray*}
\epsilon&\leq & \|\nabla f(x^{+})+g_{\psi}(x^{+})\|\\
            &\leq & \|\nabla f(x^{+})-\nabla\Phi_{x}(x^{+})\|+\|\nabla\Phi_{x}(x^{+})-\nabla \Omega_{x,M}(x^{+})\|\\
            &       & +\|\nabla\Omega_{x,M}(x^{+})+g_{\psi}(x^{+})\|\\
            &\leq &\dfrac{L_{f}}{6}\|x^{+}-x\|^{3}+\dfrac{M}{2}\|x^{+}-x\|^{3}+\dfrac{\epsilon}{2}.
\end{eqnarray*}
Thus, 
\begin{equation*}
\dfrac{\epsilon}{2}\leq\left(\dfrac{2L_{f}+6M}{12}\right)\|x^{+}-x\|^{3}=\left(\dfrac{L_{f}+3M}{6}\right)\|x^{+}-x\|^{3},
\end{equation*}
which gives
\begin{equation}
\dfrac{1}{2}\left(\dfrac{M}{L_{f}+3M}\right)\epsilon\leq\dfrac{M}{6}\|x^{+}-x\|^{3}.
\label{eq:2.19}
\end{equation}
Finally, combining (\ref{eq:2.18}) and (\ref{eq:2.19}) we conclude that (\ref{eq:2.13}) holds.
\end{proof}

\begin{remark}
\label{rem:2.0}
If $M\geq 4L_{f}$ and
\begin{equation*}
\|\nabla\Omega_{x,M}(x^{+})+g_{\psi}(x^{+})\|\leq\dfrac{\epsilon}{7}
\end{equation*}
then
\begin{equation*}
\|\nabla\Omega_{x,M}(x^{+})+g_{\psi}(x^{+})\|\leq\dfrac{1}{2}\min\left\{1,\dfrac{M}{L_{f}+3M}\right\}\epsilon.
\end{equation*}
In this case, if 
\begin{equation*}
\|\nabla f(x^{+})+g_{\psi}(x^{+})\|>\epsilon,
\end{equation*}
then it follows from Lemma \ref{lem:2.2} that $x^{+}$ satisfies (\ref{eq:2.13}).
\end{remark}

Our next lemma establishes the strong relative smoothness of $\Omega_{x,M}(\,\cdot\,)$ when $M\geq 4L_{f}$. This result will allow the use of a Bregman gradient method to compute $x^{+}$ satisfying (\ref{eq:2.13}).

\begin{lemma}
\label{lem:2.3}
Suppose that $f:\mathbb{R}^{n}\to\mathbb{R}$ satisfies (\ref{eq:2.2}) and let $M\geq 4L_{f}$. Then $\Omega_{x,M}(\,\cdot\,)$ satisfies the strong relative smoothness condition 
\begin{equation}
\dfrac{1}{2}\nabla^{2}\rho_{x}(y)\preceq\nabla^{2}\Omega_{x,M}(y)\preceq\dfrac{3}{2}\nabla^{2}\rho_{x}(y),\quad\forall x,y\in\mathbb{R}^{n},
\label{eq:2.20}
\end{equation}
for 
\begin{equation}
\rho_{x}(y)=\dfrac{1}{2}\langle\nabla^{2}f(x)(y-x),y-x\rangle+\dfrac{M}{2}d_{4}(y-x).
\label{eq:2.21}
\end{equation}
\end{lemma}

\begin{proof}
From the proof of Lemma 4 in \cite{NES1}, for $\tau^{2}=M/L_{f}$ we have
\begin{eqnarray*}
\nabla^{2}\Omega_{x,M}(y)&\succeq & \left(1-\dfrac{1}{\tau}\right)\nabla^{2}f(x)+\dfrac{M-\tau L_{f}}{2}\nabla^{2}d_{4}(y-x)\\
                                              & =  & \left(1-\dfrac{1}{\tau}\right)\left[\nabla^{2}f(x)+\dfrac{\tau M-\tau^{2}L_{f}}{(\tau-1)2}\nabla^{2}d_{4}(y-x)\right]\\
                                              & =  & \left(1-\dfrac{1}{\tau}\right)\left[\nabla^{2}f(x)+\dfrac{\tau M-M}{(\tau-1)2}\nabla^{2}d_{4}(y-x)\right]\\
                                              & = & \left(1-\dfrac{1}{\tau}\right)\left[\nabla^{2}f(x)+\dfrac{M}{2}\nabla^{2}d_{4}(y-x)\right]\\
                                             & =  &\left(1-\dfrac{1}{\tau}\right)\nabla^{2}\rho_{x}(y).
\end{eqnarray*}
Since $M\geq 4L_{f}$, it follows that $\tau\geq 2$ and so
\begin{equation}
\nabla^{2}\Omega_{x,M}(y)\succeq\dfrac{1}{2}\nabla^{2}\rho_{x}(y).
\label{eq:2.22}
\end{equation}
On the other hand, from the proof of Lemma 4 in \cite{NES1} we also have
\begin{eqnarray*}
\nabla^{2}\Omega_{x,M}(y)&\preceq & \left(\dfrac{\tau+1}{\tau-1}\right)\left[\left(1-\dfrac{1}{\tau}\right)\nabla^{2}f(x)+\dfrac{\tau(\tau-1)L_{f}}{2}\nabla^{2}d_{4}(y-x)\right]\\
                                             & = &\left(\dfrac{\tau+1}{\tau-1}\right)\left[\dfrac{\tau-1}{\tau}\nabla^{2}f(x)+\dfrac{\tau(\tau-1)L_{f}}{2}\nabla^{2}d_{4}(y-x)\right]\\
& = & \dfrac{\tau+1}{\tau}\nabla^{2}f(x)+\dfrac{(\tau+1)\tau}{2}L_{f}\nabla^{2}d_{4}(y-x)\\
&  = &\left(\dfrac{\tau+1}{\tau}\right)\left[\nabla^{2}f(x)+\dfrac{\tau^{2}L_{f}}{2}\nabla^{2}d_{4}(y-x)\right]\\
& = & \left(1+\dfrac{1}{\tau}\right)\nabla^{2}\rho_{x}(y).
\end{eqnarray*}
Since $\tau\geq 2$, it follows that
\begin{equation}
\nabla^{2}\Omega_{x,M}(y)\preceq\dfrac{3}{2}\nabla^{2}\rho_{x}(y).
\label{eq:2.23}
\end{equation}
Thus, combining (\ref{eq:2.22}) and (\ref{eq:2.23}) we get (\ref{eq:2.20}).
\end{proof}

\begin{remark}
\label{rem:2.1}
By Theorem 6  in \cite{NES0}, function $d_{4}(\,\cdot\,)$ is uniformly convex of degree 4 with parameter $1/3$, i.e., 
\begin{equation*}
d_{4}(y)\geq d_{4}(x)+\langle\nabla d_{4}(x),y-x\rangle+\dfrac{1}{4}\left(\dfrac{1}{3}\right)\|y-x\|^{4},\quad\forall x,y\in\mathbb{R}^{n}.
\end{equation*}
Consequently, function $\rho_{x}(\,\cdot\,)$ defined in (\ref{eq:2.21}) is also uniformly convex of degree $4$ with parameter $M/6$.
\end{remark}

\begin{lemma}
Suppose that $f:\mathbb{R}^{n}\to\mathbb{R}$ satisfies (\ref{eq:2.2}) and let $M\geq 4L_{f}$. Then, $\Omega_{x,M}(\,\cdot\,)$ is uniformly convex of degree 4 with parameter $M/24$.
\label{lem:2.4}
\end{lemma}

\begin{proof}
By Remark \ref{rem:2.1}, $M\geq4L_{f}$ and Lemma \ref{lem:2.3}, it follows that
\begin{eqnarray*}
\dfrac{M}{12}\|z-w\|^{4}&=&\dfrac{2\left(\frac{M}{6}\right)}{4}\|z-w\|^{4}\\
                                        &\leq & \left\la\nabla \rho_{x}(z)-\nabla\rho_{x}(w),z-w\right\ra\\
                                        & = &\left\la\int_{0}^{1}\nabla^{2}\rho_{x}\left(w+t(z-w)\right)(z-w)\,\text{d}t,z-w\right\ra\\
                                        & = &\int_{0}^{1}\left\la\nabla^{2}\rho_{x}\left(w+t(z-w)\right)(z-w),z-w\right\ra\,\text{d}t\\
                                        &\leq &\int_{0}^{1}2\left\la \nabla^{2}\Omega_{x,M}\left(w+t(z-w)\right)(z-w),z-w\right\ra\,\text{d}t\\
                                        & = &2\left\la\int_{0}^{1}\nabla^{2}\Omega_{x,M}\left(w+t(z-w)\right)(z-w)\,\text{d}t,z-w\right\ra\\
                                        & = & 2\left\la\nabla\Omega_{x,M}(z)-\nabla\Omega_{x,M}(w),z-w\right\ra.
\end{eqnarray*}
Thus, by Lemma 1 in \cite{NES2} we conclude that $\Omega_{x,M}(\,\cdot\,)$ is uniformly convex of degree 4 with parameter $M/24$.
\end{proof}

Now, given $x\in\dom\psi$ and $M>0$, consider the sublevel set
\begin{equation*}
\mathcal{L}_{M}(x)=\left\{y\in\mathbb{R}^{n}\,:\,\tilde{\Omega}_{x,M}(y)\leq\tilde{f}(x)\right\}.
\end{equation*}
The next lemma gives an upper bound to $\|y-x\|$ for $y\in\mathcal{L}_{M}(x)$.

\begin{lemma}
Suppose that $f:\mathbb{R}^{n}\to\mathbb{R}$ satisfies (\ref{eq:2.2}) and let $M\geq 4L_{f}$. Then, given $x\in\mathbb{R}^{n}$,
\begin{equation}
\|y-x\|^{3}\leq \dfrac{96\|\nabla\tilde{f}(x)\|}{M}\quad\forall y\in\mathcal{L}_{M}(x).
\label{eq:2.24}
\end{equation}
where $\nabla\tilde{f}(x)$ is an arbitrary subgradient of $\tilde{f}(\,\cdot\,)$.
\label{lem:2.5}
\end{lemma}

\begin{proof}
Let $y\in\mathcal{L}_{M}(x)$. For any $\nabla\tilde{f}(x)\in\partial\tilde{f}(x)$ we can write
\begin{equation*}
\nabla\tilde{f}(x)=\nabla f(x)+g_{\psi}(x),
\end{equation*}
where $g_{\psi}(x)\in\partial\psi(x)$. In particular, 
\begin{equation}
\psi(y)\geq\psi(x)+\la g_{\psi}(x),y-x\ra.
\label{eq:2.25}
\end{equation}
Moreover, by Lemma \ref{lem:2.4}, $\Omega_{x,M}(\,\cdot\,)$ is uniformly convex of degree 4 with parameter $M/24$. Therefore,
\begin{equation}
\Omega_{x,M}(y)\geq\Omega_{x,M}(x)+\la\nabla\Omega_{x,M}(x),y-x\ra+\dfrac{M}{96}\|y-x\|^{4}.
\label{eq:2.26}
\end{equation}
Summing up inequalities (\ref{eq:2.25}) and (\ref{eq:2.26}) and using the equality $\nabla\Omega_{x,M}(x)=\nabla f(x)$, we get
\begin{equation}
\tilde{\Omega}_{x,M}(y)\geq\tilde{\Omega}_{x,M}(x)+\la\nabla\tilde{f}(x),y-x\ra+\dfrac{M}{96}\|y-x\|^{4}.
\label{eq:2.27}
\end{equation}
By assumption, $\tilde{\Omega}_{x,M}(y)\leq\tilde{f}(x)=\tilde{\Omega}_{x,M}(x)$. Thus, subtracting $\tilde{\Omega}_{x,M}(x)$ in both sides of (\ref{eq:2.27}) it follows that
\begin{equation*}
0\geq\la\nabla\tilde{f}(x),y-x\ra+\dfrac{M}{96}\|y-x\|^{4}.
\end{equation*}
Then 
\begin{equation*}
\dfrac{M}{96}\|y-x\|^{4}\leq\la\nabla\tilde{f}(x),x-y\rangle\leq\|\nabla\tilde{f}(x)\|\|x-y\|,
\end{equation*}
which implies (\ref{eq:2.24}).
\end{proof}

From Lemma \ref{lem:2.5} we can obtain an upper bound to $\|\nabla^{2}\rho_{x}(y)\|$ for $y\in\mathcal{L}_{M}(x)$.

\begin{lemma}
\label{lem:2.6}
Suppose that $f:\mathbb{R}^{n}\to\mathbb{R}$ satisfies (\ref{eq:2.2}) and let $M\geq 4L_{f}$. Then, given $x\in\mathbb{R}^{n}$, we have 
\begin{equation}
\|\nabla^{2}\rho_{x}(y)\|\leq L_{x,M}, \quad\forall y\in\text{conv}\left(\mathcal{L}_{M}(x)\right).
\label{eq:2.28}
\end{equation}
with 
\begin{equation}
L_{x,M}\equiv \text{trace}\left(\nabla^{2}f(x)\right)+\dfrac{3M}{2}\left[\dfrac{96\|\nabla\tilde{f}(x)\|}{M}\right]^{\frac{2}{3}},
\label{eq:2.29}
\end{equation}
where $\nabla\tilde{f}(x)$ is an arbitrary subgradient of $\tilde{f}(\,\cdot\,)$.
\end{lemma}

\begin{proof}
Indeed, by (\ref{eq:2.21}) we have
\begin{equation*}
\|\nabla^{2}\rho_{x}(y)\|=\|\nabla^{2}f(x)+\dfrac{M}{2}\nabla^{2}d_{4}(y)\|\leq \|\nabla^{2}f(x)\|+\dfrac{3M}{2}\|y-x\|^{2}.
\end{equation*}
Thus, by Lemma \ref{lem:2.5} and (\ref{eq:2.29}) we have
\begin{equation*}
\|\nabla^{2}\rho_{x}(y)\|\leq L_{x,M},\quad \forall y\in\mathcal{L}_{M}(x).
\end{equation*}
Finally, since $\tilde{\Omega}_{x,M}(\,\cdot\,)$ is convex (by Lemma \ref{lem:2.3}), it follows that $\mathcal{L}_{M}(x)=\text{conv}\left(\mathcal{L}_{M}(x)\right)$.
\end{proof}

In view of Lemma \ref{lem:2.3}, Remark \ref{rem:2.1} and Lemma \ref{lem:2.6}, we can approximately solve (\ref{eq:2.12}) using the following Bregman gradient method (Algorithm B in \cite{GN1}):
\\[0.2cm]
\begin{mdframed}
\begin{equation}
y_{k+1}=\arg\min_{y\in\mathbb{R}^{n}}\left\{\langle\nabla\Omega_{x,M}(y_{k}),y-y_{k}\rangle+3\beta_{\rho_{x}}(y_{k},y)+\psi(y)\right\},\quad k\geq 0,
\label{eq:2.30}
\end{equation}
\end{mdframed}
\vspace{0.2cm}
where 
\begin{equation}
\beta_{\rho_{x}}(u,v)=\rho_{x}(v)-\rho_{x}(u)-\langle\nabla\rho_{x}(u),v-u\rangle.
\label{eq:2.31}
\end{equation}
Specifically, when $M$ is sufficiently large we have the following convergence rate.

\begin{lemma}
\label{lem:2.7}
Suppose that $f:\mathbb{R}^{n}\to\mathbb{R}$ satisfies (\ref{eq:2.2}). Given $x\in\mathbb{R}^{n}$ and $M>0$, let $\left\{y_{k}\right\}_{k\geq 0}$ be the sequence generated by method (\ref{eq:2.30}) with $y_{0}=x$. If $M\geq 4L_{f}$ then, for all $k\geq 0$,
\begin{equation}
g_{\psi}(y_{k+1})\equiv -\nabla\Omega_{x,M}(y_{k})+3\left[\nabla\rho_{x}(y_{k})-\nabla\rho_{x}(y_{k+1})\right]\in\partial \psi(y_{k+1}),
\label{eq:2.32}
\end{equation}
and
\begin{equation}
\|\nabla\Omega_{x,M}(y_{k+1})+g_{\psi}(y_{k+1})\|^{4}\leq\dfrac{3^{8} L_{x,M}^{4}\beta_{x,M}}{2M\left(\frac{6}{5}\right)^{k}},
\label{eq:2.33}
\end{equation}
where $L_{x,M}$ is given in (\ref{eq:2.29}) and 
\begin{equation}
\beta_{x,M}\equiv \dfrac{1}{2}\text{trace}\left(\nabla^{2}f(x)\right)\left[\dfrac{96\|\nabla \tilde{f}(x)\|}{M}\right]^{\frac{2}{3}}+\dfrac{M}{8}\left[\dfrac{96\|\nabla \tilde{f}(x)\|}{M}\right]^{\frac{4}{3}},
\label{eq:2.34}
\end{equation}
with $\nabla\tilde{f}(x)$ being an arbitrary subgradient of $\tilde{f}(\,\cdot\,)$. 
\end{lemma}

\begin{proof}
By Lemma A.8 in \cite{GN1}, (\ref{eq:2.32}) holds and 
\begin{equation}
\tilde{\Omega}_{x,M}(y_{k})-\tilde{\Omega}_{x,M}(y_{k+1})\geq\dfrac{\left(\frac{M}{6}\right)}{4\left(\frac{3}{2}\right)^{3}\left(3L_{x,M}\right)^{4}}\|\nabla\Omega_{x,M}(y_{k+1})+g_{\psi}(y_{k+1})\|^{4}.
\label{eq:2.35}
\end{equation}
On the other hand, by Lemma A.9 in \cite{GN1} we also have
\begin{equation}
\tilde{\Omega}_{x,M}(y_{k})-\tilde{\Omega}_{x,M}(S_{M}(x))\leq\dfrac{\left(\frac{1}{2}\right)\beta_{\rho_{x}}(x,S_{M}(x))}{\left(\frac{3}{3-\frac{1}{2}}\right)^{k}}
\label{eq:2.36}
\end{equation}
where $S_{M}(x)=\arg\min_{y\in\mathbb{R}^{n}}\tilde{\Omega}_{x,M}(y)$. Then, combining (\ref{eq:2.35}) and (\ref{eq:2.36}) we obtain
\begin{equation}
\|\nabla\Omega_{x,M}(y_{k+1})+g_{\psi}(y_{k+1})\|^{4}\leq\dfrac{3^{8} L_{x,M}^{4}\beta_{\rho_{x}}(x,S_{M}(x))}{2M\left(\frac{6}{5}\right)^{k}}.
\label{eq:2.37}
\end{equation}
Notice that $S_{M}(x)\in\mathcal{L}_{M}(x)$. Thus, it follows from (\ref{eq:2.31}), (\ref{eq:2.21}) and Lemma \ref{lem:2.5} that
\begin{eqnarray}
\beta_{\rho_{x}}(x,S_{M}(x))&=&\rho_{x}(S_{M}(x))-\rho_{x}(x)-\langle\nabla\rho_{x}(x),S_{M}(x)-x\rangle= \rho_{x}(S_{M}(x))\nonumber\\
                                          & = &\dfrac{1}{2}\langle\nabla^{2}f(x)(S_{M}(x)-x),S_{M}(x)-x\rangle+\dfrac{M}{8}\|S_{M}(x)-x\|^{4}\nonumber\\
                                          &\leq & \dfrac{1}{2}\|\nabla^{2}f(x)\|\|S_{M}(x)-x\|^{2}+\dfrac{M}{8}\|S_{M}(x)-x\|^{4}\nonumber\\
                                          &\leq & \dfrac{1}{2}\|\nabla^{2}f(x)\|\left[\dfrac{96\|\nabla \tilde{f}(x)\|}{M}\right]^{\frac{2}{3}}+\dfrac{M}{8}\left[\dfrac{96\|\nabla \tilde{f}(x)\|}{M}\right]^{\frac{4}{3}}\nonumber\\
&\leq  &\beta_{x,M}.
\label{eq:2.38}
\end{eqnarray}
Finally, combining (\ref{eq:2.37}) and (\ref{eq:2.38}) we get (\ref{eq:2.33}) holds.
\end{proof}

For any $x\in\mathbb{R}^{n}$ and $M>0$ we can apply method (\ref{eq:2.30}) with $y_{0}=x$ to try to find an approximate stationary point of $\tilde{\Omega}_{x,M}(\,\cdot\,)$ in the sense of (\ref{eq:2.13}). By Lemma \ref{lem:2.7}, if for some $k$ inequality (\ref{eq:2.33}) is not satisfied, then $M<4L_{f}$. This fact motivates the following algorithm.
\\[0.2cm]
\begin{mdframed}
\noindent\textbf{Algorithm 1}
\\[0.2cm]
\noindent\textbf{Step 0.} Given $x\in\dom\psi$, $\nabla\tilde{f}(x)\in\partial\tilde{f}(x)$ and $M,\epsilon>0$, compute
\begin{equation*}
L_{x,M}= \text{trace}\left(\nabla^{2}f(x)\right)+\dfrac{3M}{2}\left[\dfrac{96\|\nabla\tilde{f}(x)\|}{M}\right]^{\frac{2}{3}},
\end{equation*} 
and
\begin{equation*}
\beta_{x,M}= \dfrac{1}{2}\text{trace}\left(\nabla^{2}f(x)\right)\left[\dfrac{96\|\nabla \tilde{f}(x)\|}{M}\right]^{\frac{2}{3}}+\dfrac{M}{8}\left[\dfrac{96\|\nabla \tilde{f}(x)\|}{M}\right]^{\frac{4}{3}}.
\end{equation*}
Set $y_{0}=x$, $\alpha:=0$ and $k:=0$.
\\[0.2cm]
\noindent\textbf{Step 1.} Compute
\begin{equation}
y_{k+1}=\arg\min_{y\in\mathbb{R}^{n}}\left\{\langle\nabla\Omega_{x,M}(y_{k}),y-y_{k}\rangle+3\beta_{\rho_{x}}(y_{k},y)+\psi(y)\right\},
\label{eq:2.41}
\end{equation}
and 
\begin{equation*}
g_{\psi}(y_{k+1})= -\nabla\Omega_{x,M}(y_{k})+3\left[\nabla\rho_{x}(y_{k})-\nabla\rho_{x}(y_{k+1})\right].
\end{equation*}
\noindent\textbf{Step 2.} If 
\begin{equation}
\|\nabla\Omega_{x,M}(y_{k+1})+g_{\psi}(y_{k+1})\|\leq\dfrac{\epsilon}{7}\quad\text{or}\quad \|\nabla\Omega_{x,M}(y_{k+1})+g_{\psi}(y_{k+1})\|\leq\dfrac{M}{6}\|y_{k+1}-x\|^{3},
\label{eq:2.39}
\end{equation}
set $x^{+}=y_{k+1}$ and stop.
\\[0.2cm]
\noindent\textbf{Step 3.} If 
\begin{equation}
\|\nabla\Omega_{x,M}(y_{k+1})+g_{\psi}(y_{k+1})\|^{4}>\dfrac{3^{8} L_{x,M}^{4}\beta_{x,M}}{2M\left(\frac{6}{5}\right)^{k}},
\label{eq:2.40}
\end{equation}
set $x^{+}=y_{k+1}$, $\alpha:=1$, and stop.
\\[0.2cm]
\noindent\textbf{Step 4.} Set $k:=k+1$, and go to Step 1.
\end{mdframed}
\vspace{0.2cm}
Notice that the output of Algorithm 1 is a pair $(x^{+},\alpha)$. When $\alpha=1$, this means that Algorithm 1 stopped at Step 3. In this case, by Lemma \ref{lem:2.7} we have a certificate that $M<4L_{f}$. The next theorem establishes a complexity bound for Algorithm 1.
 
\begin{theorem}
\label{thm:2.1}
Suppose that $f:\mathbb{R}^{n}\to\mathbb{R}$ satisfies (\ref{eq:2.2}). Then, Algorithm 1 performs at most 
\begin{equation*}
1+\log_{1.2}\left(\dfrac{3^{8}\left(7L_{x,M}\right)^{4}\beta_{x,M}}{2M}\epsilon^{-4}\right)
\end{equation*}
iterations. Moreover, if $M\geq 4L_{f}$, then Algorithm 1 returns $\alpha=0$ and $x^{+}$ satisfying $\|\nabla f(x^{+})+g_{\psi}(x^{+})\|\leq\epsilon$ or (\ref{eq:2.13}).
\end{theorem}

\begin{proof}
Let $\left\{y_{k}\right\}_{k=0}^{\ell+1}$ be generated by Algorithm 1. Since $y_{\ell+1}$ was generated, this means that Algorithm 1 have not stopped at the $\ell$th iteration. Consequently, 
\begin{equation}
\|\nabla\Omega_{x,M}(y_{\ell})+g_{\psi}(y_{\ell})\|>\dfrac{\epsilon}{7},
\label{eq:2.42}
\end{equation}
\begin{equation}
\|\nabla\Omega_{x,M}(y_{\ell})+g_{\psi}(y_{\ell})\|>\dfrac{M}{6}\|y_{\ell}-x\|^{3},
\label{eq:2.43}
\end{equation}
and 
\begin{equation}
\|\nabla\Omega_{x,M}(y_{\ell})+g_{\psi}(y_{\ell})\|^{4}\leq\dfrac{3^{8} L_{x,M}^{4}\beta_{x,M}}{2M\left(\frac{6}{5}\right)^{\ell-1}}.
\label{eq:2.44}
\end{equation}
Assume that 
\begin{equation*}
\ell\geq 1+\log_{1.2}\left(\dfrac{3^{8}\left(7L_{x,M}\right)^{4}\beta_{x,M}}{2M}\epsilon^{-4}\right).
\end{equation*}
In this case, by (\ref{eq:2.44}) we have
\begin{equation*}
\|\nabla\Omega_{x,M}(y_{\ell})+g_{\psi}(y_{\ell})\|^{4}\leq \left(\dfrac{\epsilon}{7}\right)^{4}
\end{equation*}
and so
\begin{equation}
\|\nabla\Omega_{x,M}(y_{\ell})+g_{\psi}(y_{\ell})\|\leq\dfrac{\epsilon}{7},
\label{eq:2.45}
\end{equation}
which contradicts (\ref{eq:2.42}). Thus, we must have
\begin{equation*}
\ell<1+\log_{1.2}\left(\dfrac{3^{8}\left(7L_{x,M}\right)^{4}\beta_{x,M}}{2M}\epsilon^{-4}\right).
\end{equation*}

Now, let $(x^{+},\alpha)$ be the output of Algorithm 1, and assume that $M\geq 4L_{f}$. In this case, by Lemma \ref{lem:2.7} and Remark \ref{rem:2.0}, Algorithm 1 must stop at Step 2, which gives a pair $(x^{+},\alpha)$ with $\alpha=0$ and $x^{+}$ satisfying (\ref{eq:2.13}) or $\|\nabla f(x^{+})+g_{\psi}(x^{+})\|\leq\epsilon$. 
\end{proof}

\section{Basic Adaptive Third-Order Method}

Now that we have a suitable solver for auxiliary problems of the form (\ref{eq:2.12}) we can describe our basic adaptive third-order method to solve (\ref{eq:2.1}).
\begin{mdframed}
\noindent\textbf{Algorithm 2}
\\[0.1cm]
\noindent\textbf{Step 0.} Choose $x_{0}\in\dom\psi$, $M_{0},\epsilon>0$, and set $t:=0$.
\\[0.1cm]
\noindent\textbf{Step 1.} Find the smallest $i\geq 0$ such that $2^{i}M_{t}\geq 2M_{0}$.
\\[0.1cm]
\noindent\textbf{Step 1.1.} Apply Algorithm 1 to approximately minimize $\tilde{\Omega}_{x_{t},2^{i}M_{t}}(\,\cdot\,)$. Let $(x^{+}_{t,i},\alpha_{t,i})$ be its output, and $\nabla\tilde{f}(x^{+}_{t,i})=\nabla f(x^{+}_{t,i})+g_{\psi}(x^{+}_{t,i})$.
\\[0.1cm]
\noindent\textbf{Step 1.2.} If $\alpha_{t,i}=1$, set $i:=i+1$ and go back to Step 1.1.
\\[0.1cm]
\noindent\textbf{Step 1.3.} If $\|\nabla\tilde{f}(x_{t,i}^{+})\|\leq\epsilon$, set $x_{t+1}=x_{t,i}^{+}$ and stop.
\\[0.1cm]
\noindent\textbf{Step 1.4.} If 
\begin{equation}
\tilde{f}(x_{t})-\tilde{f}(x_{t,i}^{+})\geq\dfrac{1}{6\left(2^{i}M_{t}\right)^{\frac{1}{3}}}\|\nabla \tilde{f}(x_{t,i}^{+})\|^{\frac{4}{3}}
\label{eq:3.1}
\end{equation}
holds, set $i_{t}=i$ and go to Step 2. Otherwise, set $i:=i+1$ and go to Step 1.1.
\\[0.1cm]
\noindent\textbf{Step 2.} Set $x_{t+1}=x_{t,i}^{+}$, $M_{t+1}=2^{i_{t}-1}M_{t}$, $t:=t+1$ and go to Step 1.
\end{mdframed}
\vspace{0.2cm}
If $i$ is sufficiently large we have 
\begin{equation*}
2^{i}M_{t}\geq 4L_{f}.
\end{equation*}
Thus, in view of Lemma \ref{lem:2.1} and Theorem \ref{thm:2.1}, either $\|\nabla \tilde{f}(x_{t,i})\|\leq\epsilon$ or (\ref{eq:3.1}) holds. The next lemma gives an upper bound for $\left\{M_{t}\right\}_{t=0}^{T}$ and also a bound for the total number of executions of Algorithm 1 up to the $T$th iteration of Algorithm 2.

\begin{lemma}
\label{lem:3.1}
Suppose that $f:\mathbb{R}^{n}\to\mathbb{R}$ satisfies (\ref{eq:2.2}). Assume that $\left\{x_{t}\right\}_{t=0}^{T}$ is generated by Algorithm 2 with 
\begin{equation}
\|\nabla \tilde{f}(x^{+}_{t,i})\|>\epsilon,\quad i=0,\ldots,i_{t}\quad\text{and}\quad t=0,\ldots,T.
\label{eq:3.2}
\end{equation}
Then
\begin{equation}
M_{0}\leq M_{t}\leq\max\left\{2M_{0},4L_{f}\right\},
\label{eq:3.3}
\end{equation}
for $t=0,\ldots,T$. Moreover, the total number $E_{T}$ of executions of Algorithm 1 up to the $T$th iteration of Algorithm 2 is bounded as follows:
\begin{equation}
E_{T}\leq 2(T+1)+\log_{2}\left(\max\left\{2M_{0},4L_{f}\right\}\right)-\log_{2}(M_{0}).
\label{eq:3.4}
\end{equation}
\end{lemma} 
\begin{proof}
First let us proof that (\ref{eq:3.3}) is true for $t=0,\ldots,T$. Clearly, it is true for $t=0$. Now, given $t\geq 0$, suppose that
\begin{equation}
2^{i_{t}}M_{t}>2\max\left\{2M_{0},4L_{f}\right\}.
\label{eq:3.5}
\end{equation}
Then,
\begin{equation}
2^{i_{t}-1}M_{t}>\max\left\{2M_{0},4L_{f}\right\},
\label{eq:3.6}
\end{equation}
and, by Theorem \ref{thm:2.1} and (\ref{eq:3.2}) we have
\begin{equation}
\|\nabla\Omega_{x_{t},2^{i_{t}-1}M_{t}}(x^{+}_{t,i_{t}-1})+g_{\psi}(x^{+}_{t,i_{t}-1})\|\leq\dfrac{2^{i_{t}-1}M_{t}}{6}\|x_{t,i_{t}-1}^{+}-x_{t}\|^{3}.
\label{eq:3.7}
\end{equation}
In this case, using the convexity of $\tilde{f}(\,\cdot\,)$ it follows from (\ref{eq:3.7}) and Lemma \ref{lem:2.1} that
\begin{equation*}
\tilde{f}(x_{t})-\tilde{f}(x_{t,i_{t}-1}^{+})\geq\dfrac{1}{6\left(2^{i_{t}-1}M_{t}\right)}\|\nabla \tilde{f}(x_{t,i_{t}-1}^{+})\|^{\frac{4}{3}},
\end{equation*}
which contradicts the minimality of $i_{t}$. Thus, (\ref{eq:3.5}) cannot be true, i.e., 
\begin{equation*}
2^{i_{t}}M_{t}\leq 2\max\left\{2M_{0},4L_{f}\right\},
\end{equation*}
and so by Step 1 of Algorithm 2 we have
\begin{equation*}
M_{0}\leq M_{t+1}=2^{i_{t}-1}M_{t}\leq\max\left\{2M_{0},4L_{f}\right\},
\end{equation*}
that is, (\ref{eq:3.3}) also holds for $t\geq 1$.

Now, notice that at the $t$th iteration of Algorithm 2, Algorithm 1 is executed at most $i_{t}+1$ times. From $M_{t+1}=2^{i_{t}-1}M_{t}$ we have
\begin{equation}
i_{t}-1=\log_{2}(M_{t+1})-\log_{2}(M_{t}).
\label{eq:3.8}
\end{equation}
Then, by (\ref{eq:3.8}) and (\ref{eq:3.3}) we have
\begin{eqnarray*}
E_{T}& = &\sum_{t=0}^{T}(i_{t}+1)=\sum_{t=0}^{T}\,2+\log_{2}(M_{t+1})-\log_{2}(M_{t})\\
        & =  & 2T+\log_{2}(M_{T+1})-\log_{2}(M_{0})\\
        &\leq & 2(T+1)+\log_{2}\left(\max\left\{2M_{0},4L_{f}\right\}\right)-\log_{2}(M_{0}).
\end{eqnarray*}
Thus, (\ref{eq:3.4}) is also true.
\end{proof}

In what follows we will assume that 
\begin{equation}
\max_{x\in\mathcal{L}(x_{0})}\|x-x^{*}\|\leq R_{0}\in (0,+\infty),
\label{eq:3.9}
\end{equation}
where 
\begin{equation*}
\mathcal{L}(x_{0})=\left\{x\in\mathbb{R}^{n}\,:\,\tilde{f}(x)\leq \tilde{f}(x_{0})\right\}.
\end{equation*}
Using (\ref{eq:3.9}) we can also bound the number of inner iterations at each execution of Algorithm 1 in Algorithm 2. 

\begin{lemma}
\label{lem:3.2}
Suppose that $f:\mathbb{R}^{n}\to\mathbb{R}$ satisfies (\ref{eq:2.2}) and that (\ref{eq:3.9}) holds. Assume that $\left\{x_{t}\right\}_{t=0}^{T}$ is generated by Algorithm 2 with
\begin{equation}
\|\nabla \tilde{f}(x^{+}_{t,i})\|>\epsilon,\quad i=0,\ldots,i_{t}\quad\text{and}\quad t=0,\ldots,T.
\label{eq:3.10}
\end{equation}
Then, for each $t\in\left\{0,\ldots,T\right\}$ and each $i\in\left\{0,\ldots,i_{t}\right\}$, the computation of $x^{+}_{t,i}$ requires at most
\begin{equation*}
\mathcal{O}\left(\log_{1.2}\left(\dfrac{3^{8}\left(7L_{0}\right)^{4}\beta_{0}}{4M_{0}}\epsilon^{-4}\right)\right)
\end{equation*}
iterations of Algorithm 1, where
\begin{equation}
L_{0}=4nL_{f}R_{0}^{2}+2\text{trace}\left(\nabla^{2}f(x_{0})\right)+3\max\left\{2M_{0},4L_{f}\right\}\left(\dfrac{48}{M_{0}}\right)^{\frac{2}{3}}G_{0}^{\frac{2}{3}},
\label{eq:3.11}
\end{equation}
\begin{equation}
\beta_{0}=\left(2nL_{f}R_{0}^{2}+\text{trace}\left(\nabla^{2}f(x_{0})\right)\right)\left(\dfrac{48}{M_{0}}\right)^{\frac{2}{3}}G_{0}^{\frac{2}{3}}+\dfrac{\max\left\{2M_{0},4L_{f}\right\}}{4}\left(\dfrac{48}{M_{0}}\right)^{\frac{4}{3}}G_{0}^{\frac{4}{3}},
\label{eq:3.12}
\end{equation}
and
\begin{equation}
G_{0}=\max\left\{\|\nabla\tilde{f}(x_{0})\|,6^{\frac{3}{4}}\max\left\{4M_{0},8L_{f}\right\}^{\frac{1}{4}}(\tilde{f}(x_{0})-\tilde{f}(x^{*}))^{\frac{3}{4}}\right\}.
\label{eq:3.12extra}
\end{equation}
\end{lemma}

\begin{proof}
By Theorem \ref{thm:2.1}, the computation of $x^{+}_{t,i}$ requires at most
\begin{equation*}
\mathcal{O}\left(\log_{1.2}\left(\dfrac{3^{8}\left(7 L_{x_{t},2^{i}M_{t}}\right)^{4}\beta_{x_{t},2^{i}M_{t}}}{2(2^{i}M_{t})}\epsilon^{-4}\right)\right),
\end{equation*}
where 
\begin{equation}
L_{x_{t},2^{i}M_{t}}=\text{trace}\left(\nabla^{2}f(x_{t})\right)+\dfrac{3(2^{i}M_{t})}{2}\left(\dfrac{96}{2^{i}M_{t}}\right)^{\frac{2}{3}}\|\nabla\tilde{f}(x_{t})\|^{\frac{2}{3}}
\label{eq:3.13}
\end{equation}
and
\begin{equation}
\beta_{x_{t},2^{i}M_{t}}=\dfrac{1}{2}\text{trace}\left(\nabla^{2}f(x_{t})\right)\left(\dfrac{96}{2^{i}M_{t}}\right)^{\frac{2}{3}}\|\nabla\tilde{f}(x_{t})\|^{\frac{2}{3}}+\dfrac{2^{i}M_{t}}{8}\left(\dfrac{96}{2^{i}M_{t}}\right)^{\frac{4}{3}}\|\nabla\tilde{f}(x_{t})\|^{\frac{4}{3}}.
\label{eq:3.14}
\end{equation}
Regarding $2^{i}M_{t}$, it follows from Step 1 and Lemma \ref{lem:3.1} that
\begin{equation}
2M_{0}\leq 2^{i}M_{t}\leq 2^{i_{t}}M_{t}=2M_{t+1}\leq 2\max\left\{2M_{0},4L_{f}\right\}.
\label{eq:3.15}
\end{equation}
Thus, the computation of $x^{+}_{t,i}$ requires at most
\begin{equation*}
\mathcal{O}\left(\log_{1.2}\left(\dfrac{3^{8}\left(7L_{x_{t},2^{i}M_{t}}\right)^{4}\beta_{x_{t},2^{i}M_{t}}}{4M_{0}}\epsilon^{-4}\right)\right)
\end{equation*}
It remains to show that $L_{x_{t},2^{i}M_{t}}$ and $\beta_{x_{t},2^{i}M_{t}}$ are bounded from above by $L_{0}$ and $\beta_{0}$, respectively. In view of (\ref{eq:2.6}) we have
\begin{eqnarray*}
\nabla^{2}f(x_{t})&=&\nabla^{2}f(x_{t})-\nabla\Phi_{x_{0}}(x_{t})+\nabla^{2}\Phi_{x_{0}}(x_{t})\\
                               &\preceq & \dfrac{L_{f}}{2}\|x_{t}-x_{0}\|^{2}I+\nabla^{2}f(x_{0})+D^{3}f(x_{0})[x_{t}-x_{0}],
\end{eqnarray*}
and so
\begin{equation}
\text{trace}\left(\nabla^{2}f(x_{t})\right)\leq\dfrac{nL_{f}}{2}\|x_{t}-x_{0}\|^{2}+\text{trace}\left(\nabla^{2}f(x_{0})\right)+\text{trace}\left(D^{3}f(x_{0})[x_{t}-x_{0}]\right).
\label{eq:3.16}
\end{equation}
On the other hand, by Lemma 3 in \cite{NES1} we also have
\begin{equation*}
D^{3}f(x_{0})[x_{t}-x_{0}]\preceq\nabla^{2}f(x_{0})+\dfrac{L_{f}}{2}\|x_{t}-x_{0}\|^{2}I,
\end{equation*}
which gives
\begin{equation}
\text{trace}\left(D^{3}f(x_{0})[x_{t}-x_{0}]\right)\leq\text{trace}\left(\nabla^{2}f(x_{0})\right)+\dfrac{n L_{f}}{2}\|x_{t}-x_{0}\|^{2}.
\label{eq:3.17}
\end{equation}
Thus, combining (\ref{eq:3.16}) and (\ref{eq:3.17}) we obtain
\begin{eqnarray*}
\text{trace}\left(\nabla^{2}f(x_{t})\right)&\leq & n L_{f}\|x_{t}-x_{0}\|^{2}+2\text{trace}\left(\nabla^{2}f(x_{0})\right)\\
                                                                    &\leq & n L_{f}\left(\|x_{t}-x^{*}\|+\|x^{*}-x_{0}\|\right)^{2}+2\text{trace}\left(\nabla^{2}f(x_{0})\right)\\
                                                                    &\leq & 2n L_{f}\left(\|x_{t}-x^{*}\|^{2}+\|x_{0}-x^{*}\|^{2}\right)+2\text{trace}\left(\nabla^{2}f(x_{0})\right).
\end{eqnarray*}
Since $x_{t}\in\mathcal{L}(x_{0})$, it follows from (\ref{eq:3.9}) that
\begin{equation}
\text{trace}\left(\nabla^{2}f(x_{t})\right)\leq 4n L_{f}R_{0}^{2}+2\text{trace}\left(\nabla^{2}f(x_{0})\right).
\label{eq:3.18}
\end{equation}
Now, suppose that $t\geq 1$. By (\ref{eq:3.1}) we have
\begin{equation*}
\tilde{f}(x_{0})-\tilde{f}(x^{*})\geq \tilde{f}(x_{t-1})-\tilde{f}(x_{t})\geq\dfrac{1}{6(2M_{t+1})^{\frac{1}{3}}}\|\nabla\tilde{f}(x_{t})\|^{\frac{4}{3}}.
\end{equation*}
Thus, by Lemma \ref{lem:3.1}, we get
\begin{equation}
\|\nabla\tilde{f}(x_{t})\|^{\frac{4}{3}}\leq 6\max\left\{4M_{0},8L_{f}\right\}^{\frac{1}{3}}(\tilde{f}(x_{0})-\tilde{f}(x^{*}))
\label{eq:3.19}
\end{equation}
and so
\begin{equation}
\|\nabla\tilde{f}(x_{t})\|\leq 6^{\frac{3}{4}}\max\left\{4M_{0},8L_{f}\right\}^{\frac{1}{4}}(\tilde{f}(x_{0})-\tilde{f}(x^{*}))^{\frac{3}{4}}\leq G_{0},
\label{eq:3.20}
\end{equation}
where $G_{0}$ is given in (\ref{eq:3.12extra}). Finally, it follows from (\ref{eq:3.13}), (\ref{eq:3.15}), (\ref{eq:3.18}) and (\ref{eq:3.20}) that $L_{x_{t},2^{i}M_{t}}\leq L_{0}$, while (\ref{eq:3.14}), (\ref{eq:3.15}), (\ref{eq:3.18}) and (\ref{eq:3.19}) imply that $\beta_{x_{t},2^{i}M_{t}}\leq \beta_{0}$.
\end{proof}

\begin{remark}
From (\ref{eq:3.11}) and (\ref{eq:3.12}), it follows that $L_{0}$ and $\beta_{0}$ depend linearly on $n$, the problem dimension. As we can see in (\ref{eq:3.16}), this dependence is due to the use of $\text{trace}\left(\nabla^{2}f(x)\right)$ in the definitions of $L_{x,M}$ and $\beta_{x,M}$ at Step 0 of Algorithm 1. However, looking at the proofs of Lemmas \ref{lem:2.6} and \ref{lem:2.7}, it is clear that we could use $\|\nabla^{2}f(x)\|$ instead of $\text{trace}\left(\nabla^{2}f(x)\right)$, which would yields $L_{0}$ and $\beta_{0}$ independent of $n$. In Algorithm 1, we choose to use $\text{trace}\left(\nabla^{2}f(x)\right)$ because it can be easily computed. Since $L_{0}$ and $\beta_{0}$ appear inside the logarithm in our iteration complexity bound for Algorithm 1, the dependence on $n$ is almost negligible.

\end{remark}

\begin{theorem}
\label{thm:3.1}
Suppose that $f:\mathbb{R}^{n}\to\mathbb{R}$ satisfies (\ref{eq:2.2}) and that (\ref{eq:3.9}) holds. Given $\epsilon>0$, assume that $\left\{x_{t}\right\}_{t=0}^{T}$ is generated by Algorithm 2 with 
\begin{equation}
\min\left\{\|\nabla \tilde{f}(x_{t,i}^{+})\|,\tilde{f}(x_{t,i}^{+})-\tilde{f}(x^{*})\right\}\geq\epsilon,\quad i=0,\ldots,i_{t}\quad\text{and}\quad t=0,\ldots,T.
\label{eq:3.21}
\end{equation}
Denote by $m$ the first iteration number such that
\begin{equation}
\tilde{f}(x_{m})-\tilde{f}(x^{*})\leq 6^{3}\left(4\max\left\{2M_{0},4L_{f}\right\}R_{0}^{4}\right),
\label{eq:3.22}
\end{equation}
and assume that $T>m$. Then, 
\begin{equation}
m\leq 1+\dfrac{1}{\log_{2}\left(\frac{4}{3}\right)}\log_{2}\log_{2}\left(\dfrac{\tilde{f}(x_{0})-\tilde{f}(x^{*})}{6^{3}\left(\max\left\{4M_{0},8L_{f}\right\}R_{0}^{4}\right)}\right)
\label{eq:3.23}
\end{equation}
and, for all $k$, $m\leq k\leq T$ we have
\begin{equation}
\tilde{f}(x_{k})-\tilde{f}(x^{*})\leq\dfrac{(42)^{3}\max\left\{4M_{0},8L_{f}\right\}R_{0}^{4}}{(k-m)^{3}}.
\label{eq:3.24}
\end{equation}
Consequently,
\begin{equation}
T\leq m+42\max\left\{4M_{0},8L_{f}\right\}^{\frac{1}{3}}R_{0}^{\frac{4}{3}}\epsilon^{-\frac{1}{3}}.
\label{eq:3.25}
\end{equation}
\end{theorem}

\begin{proof}
By (\ref{eq:3.21}) and Steps 1 and 2 of Algorithm 2, we have
\begin{equation}
\tilde{f}(x_{k})-\tilde{f}(x_{k+1})\geq\dfrac{1}{6(2M_{k+1})^{\frac{1}{3}}}\|\nabla f(x_{k+1})\|^{\frac{4}{3}}\quad\text{for}\,\,k=0,\ldots,T-1.
\label{eq:3.26}
\end{equation}
Since $x_{k+1}\in\mathcal{L}(x_{0})$ for all $k\in\left\{1,\ldots,T-1\right\}$, using the convexity of $\tilde{f}(\,\cdot\,)$, (\ref{eq:3.9}) and (\ref{eq:3.3}) in (\ref{eq:3.26}) we get
\begin{eqnarray}
\tilde{f}(x_{k})-\tilde{f}(x_{k+1})&\geq & \dfrac{1}{6\left(2M_{k+1}R_{0}^{4}\right)^{\frac{1}{3}}}\left(\|\nabla \tilde{f}(x_{k+1})\|\|x_{k+1}-x^{*}\|\right)^{\frac{4}{3}}\nonumber\\
&\leq &\dfrac{1}{6\left(2M_{k+1}R_{0}^{4}\right)^{\frac{1}{3}}}\left(\tilde{f}(x_{k})-\tilde{f}(x^{*})\right)^{\frac{4}{3}}\nonumber\\
&\geq & \dfrac{1}{6\left(2\max\left\{2M_{0},4L_{f}\right\}R_{0}^{4}\right)^{\frac{1}{3}}}\left(\tilde{f}(x_{k})-\tilde{f}(x^{*})\right)^{\frac{4}{3}},\quad k=0,\ldots,T-1.\nonumber\\
\label{eq:3.27}
\end{eqnarray}
Denote 
\begin{equation}
\delta_{k}=\dfrac{\tilde{f}(x_{k})-\tilde{f}(x^{*})}{6^{3}\left(2\max\left\{2M_{0},4L_{f}\right\}R_{0}^{4}\right)}
\label{eq:3.28}
\end{equation}
Then, it follows from (\ref{eq:3.27}) that
\begin{equation}
\delta_{k}-\delta_{k+1}\geq\delta_{k+1}^{\frac{4}{3}},\quad k=0,\ldots,T-1.
\label{eq:3.29}
\end{equation}
Thus, $\left\{\delta_{k}\right\}_{k=0}^{T-1}$ satisfies condition (1.1) of Lemma 1.1 in \cite{GN0} with power equal to $4/3$. By definition, $m$ is the first iteration number such that (\ref{eq:3.11}) holds. This means that $m$ is the first iteration number for which $\delta_{m}\leq 2$. If $m>0$, then $\delta_{0}>2$. Then, by inequality (1.2) of Lemma 1.1 in \cite{GN1} we obtain
\begin{equation*}
\ln(2)\leq \ln(\delta_{m-1})\leq\left(\frac{3}{4}\right)^{m-1}\ln(\delta_{0}).
\end{equation*}
Rearranging the terms in this inequality and then taking the logarithm in both sides, we conclude that (\ref{eq:3.23}) holds. Using $\delta_{m}\leq 2$, it follows from the inequality (1.3) of Lemma 1.1 in \cite{GN0} that
\begin{equation*}
\delta_{k}\leq\left[\dfrac{1+\delta_{m}^{\frac{1}{3}}}{\frac{1}{3}(k-m)}\right]^{3}\leq \left(\dfrac{7}{k-m}\right)^{3},\quad m\leq k\leq T.
\label{eq:3.29extra}
\end{equation*}
Now, replacing (\ref{eq:3.28}) in (\ref{eq:3.29extra}) we get
\begin{equation*}
\dfrac{\tilde{f}(x_{k})-\tilde{f}(x^{*})}{6^{3}\left(2\max\left\{2M_{0},4L_{f}\right\}R_{0}^{4}\right)}\leq\dfrac{7^{3}}{(k-m)^{3}}\quad m\leq k\leq T,
\end{equation*}
which implies (\ref{eq:3.24}). 

Finally, by (\ref{eq:3.21}) and (\ref{eq:3.24}) we have the inequality
\begin{equation*}
\epsilon< \tilde{f}(x_{k})-\tilde{f}(x^{*})\leq\dfrac{(42)^{3}\max\left\{4M_{0},8L_{f}\right\}R_{0}^{4}}{(k-m)^{3}},\quad m\leq k\leq T.
\end{equation*}
from which we get the bound (\ref{eq:3.25}).
\end{proof}

Combining (\ref{eq:3.4}) and (\ref{eq:3.25}) we have the following worst-case complexity result.

\begin{corollary}
\label{cor:3.1}
Under the same assumptions of Theorem \ref{thm:3.1}, let $x_{T+1}=x^{+}_{T,i_{T}}$ be the first trial point generated by Algorithm 2 such that
\begin{equation}
\min\left\{\|\nabla \tilde{f}(x_{T+1})\|,\tilde{f}(x_{T+1})-\tilde{f}(x^{*})\right\}\leq\epsilon.
\label{eq:3.30}
\end{equation}
Then, the total number of executions of Algorithm 1 is bounded as follows
\begin{equation}
E_{T}\leq 2(m+1)+84\max\left\{4M_{0},8L_{f}\right\}^{\frac{1}{3}}R_{0}^{\frac{4}{3}}\epsilon^{-\frac{1}{3}}+\log_{2}\left(\max\left\{2M_{0},4L_{f}\right\}\right)-\log_{2}(M_{0}),
\label{eq:3.31}
\end{equation}
where $m$ is bounded as in (\ref{eq:3.23}).
\end{corollary}

\begin{remark}
Let us refer to the iterations of Algorithm 1 by \textit{inner iterations}. In view of Corollary \ref{cor:3.1} and Lemma \ref{lem:3.2}, Algorithm 2 needs at most $\mathcal{O}\left(|\log_{1.2}(\epsilon)|\epsilon^{-\frac{1}{3}}\right)$ inner iterations to generate $\bar{x}$ such that
\begin{equation*}
\min\left\{\|\nabla\tilde{f}(\bar{x})\|,\tilde{f}(\bar{x})-\tilde{f}(x^{*})\right\}\leq\epsilon.
\end{equation*}
In fact, with a direct adaptation of the proof of Theorem 3.3 in \cite{GN5}, it can be shown that an inner iteration complexity bound of the same order holds for the generation of $\bar{x}$ such that $\|\nabla\tilde{f}(\bar{x})\|\leq\epsilon$.
\end{remark}

\section{Accelerated Adaptive Third-Order Method}

The adaptive procedure used in Algorithm 2 can also be incorporated in accelerated third-order methods (see, e.g., \cite{GN01}). The next algorithm is an adaptive accelerated method based on Algorithm 1. 
\\[0.2cm]
\begin{mdframed}
\noindent\textbf{Algorithm 3}
\\[0.2cm]
\noindent\textbf{Step 0.} Choose $x_{0}\in\dom\psi$, $M_{0},\epsilon>0$, and define function $\varphi_{0}(x)=\frac{1}{4}\|x-x_{0}\|^{4}$ for all $x\in\dom\psi$. Set $v_{0}=x_{0}$, $A_{0}=0$, and $t:=0$.
\\[0.2cm]
\noindent\textbf{Step 1.} Find the smallest $i\geq 0$ such that $2^{i}M_{t}\geq 2M_{0}$.
\\[0.2cm]
\noindent\textbf{Step 1.1.} Compute the coefficient $a_{t,i}>0$ by solving the equation
\begin{equation}
a_{t,i}^{4}=\dfrac{4^{2}\left(A_{t}+a_{t,i}\right)^{3}}{18^{3}\left(2^{i}M_{t}\right)}
\label{eq:4.1}
\end{equation}
\noindent\textbf{Step 1.2.} Set $\gamma_{t,i}=\frac{a_{t,i}}{A_{t}+a_{t,i}}$ and compute the vector
\begin{equation}
z_{t,i}=(1-\gamma_{t,i})x_{t}+\gamma_{t,i}v_{t}.
\label{eq:4.2}
\end{equation}
\noindent\textbf{Step 1.3.} Apply Algorithm 1 to approximately minimize $\tilde{\Omega}_{z_{t,i},2^{i}M_{t}}(\,\cdot\,)$. Let $(x_{t,i}^{+},\alpha_{t,i})$ be its output, and $\nabla\tilde{f}(x^{+}_{t,i})=\nabla f(x^{+}_{t,i})+g_{\psi}(x^{+}_{t,i})$.
\\[0.2cm]
\noindent\textbf{Step 1.4.} If $\alpha_{t,i}=1$, set $i:=i+1$ and go back to Step 1.1.
\\[0.2cm]
\noindent\textbf{Step 1.5.} If $\|\nabla \tilde{f}(x_{t,i}^{+})\|\leq\epsilon$, set $x_{t+1}=x_{t,i}^{+}$, and stop.
\\[0.2cm]
\noindent\textbf{Step 1.6.} If
\begin{equation}
\langle\nabla\tilde{f}(x^{+}_{t,i}),z_{t,i}-x_{t,i}^{+}\rangle\geq\dfrac{1}{6(2^{i}M_{t})^{\frac{1}{3}}}\|\nabla\tilde{f}(x^{+}_{t,i})\|^{\frac{4}{3}}
\label{eq:4.3}
\end{equation}
holds, set $i_{t}=i$ and go to Step 2. Otherwise, set $i:=i+1$ and go to Step 1.1.
\\[0.2cm]
\noindent\textbf{Step 2.} Set $x_{t+1}=x^{+}_{t,i_{t}}$, $z_{t}=z_{t,i_{t}}$, $a_{t}=a_{t,i_{t}}$ and $\gamma_{t}=\gamma_{t,i_{t}}$. Define $A_{t+1}=A_{t}+a_{t}$ and $M_{t+1}=2^{i_{t}-1}M_{t}$.
\\[0.2cm]
\noindent\textbf{Step 3.} Define
\begin{equation}
\varphi_{t+1}(x)=\varphi_{t}(x)+a_{t}\left[f(x_{t+1})+\langle\nabla f(x_{t+1}),x-x_{t+1}\rangle+\psi(x)\right],\quad\forall x\in\dom\psi,
\label{eq:4.4}
\end{equation}
and compute
\begin{equation}
v_{t+1}=\arg\min_{x\in\dom\psi}\varphi_{t+1}(x).
\label{eq:4.5}
\end{equation}
\noindent\textbf{Step 4.} Set $t:=t+1$, and go back to Step 1.
\end{mdframed}
\vspace{0.2cm}

\begin{lemma}
\label{lem:4.1}
Suppose that $f:\mathbb{R}^{n}\to\mathbb{R}$ satisfies (\ref{eq:2.2}). Assume that $\left\{x_{t}\right\}_{t=0}^{T}$ is generated by Algorithm 3 with 
\begin{equation*}
\|\nabla \tilde{f}(x^{+}_{t,i})\|>\epsilon,\quad i=0,\ldots,i_{t}\quad\text{and}\quad t=0,\ldots,T.
\end{equation*}
Then
\begin{equation}
M_{0}\leq M_{t}\leq\max\left\{2M_{0},4L_{f}\right\},
\label{eq:4.6}
\end{equation}
for $t=0,\ldots,T$. Moreover, the total number $E_{T}$ of executions of Algorithm 1 up to the $T$th iteration of Algorithm 3 is bounded as follows:
\begin{equation}
E_{T}\leq 2(T+1)+\log_{2}\left(\max\left\{2M_{0},4L_{f}\right\}\right)-\log_{2}(M_{0}).
\label{eq:4.7}
\end{equation}
\end{lemma} 

\begin{proof}
By Step 1 of Algorithm 3 and Theorem \ref{thm:2.1}, the statements follow as in the proof of Lemma \ref{lem:3.1}.
\end{proof}

\begin{lemma}
\label{lem:4.2}
Suppose that $f:\mathbb{R}^{n}\to\mathbb{R}$ satisfies (\ref{eq:2.2}), and assume that $\left\{x_{t}\right\}_{t=0}^{T}$ is generated by Algorithm 2. Then, 
\begin{equation}
A_{t}\tilde{f}(x_{t})\leq\varphi_{t}^{*}\equiv\min_{x\in\dom\psi}\varphi_{t}(x),
\label{eq:4.8}
\end{equation}
for all $t\in\left\{0,\ldots,T\right\}$.
\end{lemma}

\begin{proof}
By the equality $A_{0}=0$ and the definition of $\varphi_{0}(\,\cdot\,)$, we have
\begin{equation*}
A_{0}\tilde{f}(x_{0})=0=\min_{x\in\dom\psi}\varphi_{0}(x),
\end{equation*}
that is, (\ref{eq:4.8}) is true for $t=0$. Assume that (\ref{eq:4.8}) is true for some $t\in\left\{0,\ldots,T-1\right\}$. Since $\varphi_{t}(\,\cdot\,)$ is uniformly convex of degree 4 with parameter $1/4$, by the induction assumption we have
\begin{equation}
\varphi_{t}(x)\geq\varphi_{t}^{*}+\dfrac{1}{16}\|x-v_{t}\|^{4}\geq A_{t}\tilde{f}(x_{t})+\dfrac{1}{16}\|x-v_{t}\|^{4}.
\label{eq:4.9}
\end{equation}
Then, combining (\ref{eq:4.4}) and (\ref{eq:4.9}) we get
\begin{eqnarray}
\varphi_{t+1}^{*}&=&\min_{x\in\dom\psi}\left\{\varphi_{t}(x)+a_{t}\left[f(x_{t+1})+\langle\nabla f(x_{t+1}),x-x_{t+1}\rangle+\psi(x)\right]\right\}\nonumber\\
                         &\geq &\min_{x\in\dom\psi}\left\{A_{t}\tilde{f}(x_{t})+\dfrac{1}{16}\|x-v_{t}\|^{4}+a_{t}\left[f(x_{t+1})+\langle\nabla f(x_{t+1}),x-x_{t+1}\rangle+\psi(x)\right]\right\}\nonumber\\
                         &        &
\label{eq:4.10}
\end{eqnarray}
Using the convexity of $\tilde{f}(\,\cdot\,)$, we have
\begin{equation}
A_{t}\tilde{f}(x_{t})\geq A_{t}\tilde{f}(x_{t+1})+\langle\nabla\tilde{f}(x_{t+1}),A_{t}x_{t}-A_{t}x_{t+1}\rangle.
\label{eq:4.11}
\end{equation}
Moreover, since $g_{\psi}(x_{t+1})\in\partial\psi(x_{t+1})$ we also have
\begin{equation*}
\psi(x)\geq\psi(x_{t+1})+\langle g_{\psi}(x_{t+1}),x-x_{t+1}\rangle,
\end{equation*}
which gives
\begin{equation}
a_{t}\left[f(x_{t+1})+\langle\nabla f(x_{t+1}),x-x_{t+1}\rangle+\psi(x)\right]\geq a_{t}\left[\tilde{f}(x_{t+1})+\langle\nabla\tilde{f}(x_{t+1}),x-x_{t+1}\rangle\right].
\label{eq:4.12}
\end{equation}
Thus, using (\ref{eq:4.11}), (\ref{eq:4.12}) and $A_{t+1}=A_{t}+a_{t}$ in (\ref{eq:4.10}), it follows that
\begin{equation}
\varphi_{t+1}^{*}\geq\min_{x\in\dom\psi}\left\{A_{t+1}\tilde{f}(x_{t+1})+\langle\nabla\tilde{f}(x_{t+1}),A_{t}x_{t}-A_{t}x_{t+1}+a_{t}x-a_{t}x_{t+1}\rangle+\dfrac{1}{16}\|x-v_{t}\|^{4}\right\}.
\label{eq:4.13}
\end{equation}
By (\ref{eq:4.2}) and the definitions of $\gamma_{t}$ and $A_{t+1}$, we have
\begin{equation*}
z_{t}=(1-\gamma_{t})x_{t}+\gamma_{t}v_{t}=\left(1-\dfrac{a_{t}}{A_{t}+a_{t}}\right)x_{t}+\dfrac{a_{t}}{A_{t}+a_{t}}v_{t}=\dfrac{A_{t}}{A_{t+1}}x_{t}+\dfrac{a_{t}}{A_{t+1}}v_{t}.
\end{equation*}
Consequently
\begin{equation*}
A_{t+1}z_{t}=A_{t}x_{t}+a_{t}v_{t},
\end{equation*}
which gives
\begin{equation*}
A_{t}x_{t}=A_{t+1}z_{t}-a_{t}v_{t}.
\end{equation*}
Thus,
\begin{eqnarray}
A_{t}x_{t}-A_{t}x_{t+1}+a_{t}x-a_{t}x_{t+1}&=&A_{t}x_{t}-A_{t+1}x_{t+1}+a_{t}x\nonumber\\
                                                                              &=&A_{t+1}z_{t}-a_{t}v_{t}-A_{t+1}x_{t+1}+a_{t}x\nonumber\\
                                                                              &=& A_{t+1}(z_{t}-x_{t+1})+a_{t}(x-v_{t}).
\label{eq:4.14}
\end{eqnarray}
Now, combining (\ref{eq:4.13}), (\ref{eq:4.14}) and (\ref{eq:4.3}), it follows that
\begin{eqnarray}
\varphi_{t+1}^{*}&\geq & A_{t+1}\tilde{f}(x_{t+1})+\min_{x\in\dom\psi}\left\{A_{t+1}\langle\nabla\tilde{f}(x_{t+1}),z_{t}-x_{t+1}\rangle+a_{t}\langle\nabla\tilde{f}(x_{t+1}),x-v_{t}\rangle+\dfrac{1}{16}\|x-v_{t}\|^{4}\right\}\nonumber\\
                        &\geq & A_{t+1}\tilde{f}(x_{t+1})+\min_{x\in\dom\psi}\left\{\dfrac{A_{t+1}}{6\left(2^{i_{t}}M_{t}\right)^{\frac{1}{3}}}\|\nabla\tilde{f}(x_{t+1})\|^{\frac{4}{3}}+a_{t}\langle\nabla\tilde{f}(x_{t+1}),x-v_{t}\rangle+\dfrac{1}{16}\|x-v_{t}\|^{4}\right\}\nonumber\\
                        &        &
\label{eq:4.15}
\end{eqnarray}
By Lemma 2 in \cite{NES2} with $p=4$, $s=a_{t}\nabla\tilde{f}(x_{t+1})$ and $\sigma=1/4$, we have
\begin{eqnarray*}
a_{t}\langle\nabla\tilde{f}(x_{t+1}),x-v_{t}\rangle+\dfrac{1}{16}\|x-v_{t}\|^{4}&\geq & -\left(\dfrac{3}{4}\right)4^{\frac{1}{3}}a_{t}^{\frac{4}{3}}\|\nabla\tilde{f}(x_{t+1})\|^{\frac{4}{3}}\nonumber\\
                                                                                                                                     &  =    &-\left(\dfrac{3}{4^{\frac{2}{3}}}\right)a_{t}^{\frac{4}{3}}\|\nabla\tilde{f}(x_{t+1})\|^{\frac{4}{3}}.
\end{eqnarray*}
Thus, 
\begin{eqnarray}
\quad\quad & & \min_{x\in\dom\psi}\left\{\dfrac{A_{t+1}}{6\left(2^{i_{t}}M_{t}\right)^{\frac{1}{3}}}\|\nabla\tilde{f}(x_{t+1})\|^{\frac{4}{3}}+a_{t}\langle\nabla\tilde{f}(x_{t+1}),x-v_{t}\rangle+\dfrac{1}{16}\|x-v_{t}\|^{4}\right\}\nonumber\\
\quad\quad&\geq & \left[\dfrac{A_{t+1}}{6\left(2^{i_{t}}M_{t}\right)^{\frac{1}{3}}}-4^{-\frac{2}{3}}3a_{t}^{\frac{4}{3}}\right]\|\nabla\tilde{f}(x_{t+1})\|^{\frac{4}{3}}\nonumber\\
\quad\quad &\geq & \left[\dfrac{(A_{t}+a_{t})}{6\left(2^{i_{t}}M_{t}\right)^{\frac{1}{3}}}-4^{-\frac{2}{3}}3a_{t}^{\frac{4}{3}}\right]\|\nabla\tilde{f}(x_{t+1})\|^{\frac{4}{3}}\nonumber\\
\quad\quad & =    & 0,
\label{eq:4.16}
\end{eqnarray}
where the last equality follows from (\ref{eq:4.1}). Finally, combining (\ref{eq:4.15}) and (\ref{eq:4.16}), we conclude that
\begin{equation*}
\varphi_{t+1}^{*}\geq A_{t+1}\tilde{f}(x_{t+1}),
\end{equation*}
that is, (\ref{eq:4.8}) is true for $t+1$. 
\end{proof}

\begin{lemma}
Let $\left\{x_{t}\right\}_{t=0}^{T}$ be generated by Algorithm 3. Then, for all $t\in\left\{0,\ldots,T\right\}$, 
\begin{equation*}
\varphi_{t}(x)\leq A_{t}\tilde{f}(x)+\dfrac{1}{4}\|x-x_{0}\|^{4},\quad\forall x\in\mathbb{R}^{n}.
\end{equation*}
\label{lem:4.3}
\end{lemma}

\begin{proof}
It follows by induction as in the proof of Lemma 3.2 in \cite{GN02}.
\end{proof}

Now, combining Lemmas \ref{lem:4.1}, \ref{lem:4.2} and \ref{lem:4.3} we can establish the rate of convergence of Algorithm 3.

\begin{theorem}
\label{thm:4.1}
Suppose that $f:\mathbb{R}^{n}\to\mathbb{R}$ satisfies (\ref{eq:2.2}). Assume that $\left\{x_{t}\right\}_{t=0}^{T}$, with $T\geq 2$, is generated by Algorithm 3 with 
\begin{equation}
\min\left\{\|\nabla \tilde{f}(x^{+}_{t,i})\|,\tilde{f}(x_{t})-\tilde{f}(x^{*})\right\}>\epsilon,\quad i=0,\ldots,i_{t}\quad\text{and}\quad t=0,\ldots,T.
\label{eq:4.17}
\end{equation}
Then
\begin{equation}
\tilde{f}(x_{t})-\tilde{f}(x^{*})\leq\dfrac{18^{3}\max\left\{4M_{0},8L_{f}\right\}\|x_{0}-x^{*}\|^{4}}{2^{\frac{13}{4}}(t-1)^{4}},\quad \forall t\in\left\{2,\ldots,T\right\}.
\label{eq:4.18}
\end{equation}
Consequently,
\begin{equation}
T\leq 1+\left(\frac{1}{2}\right)^{\frac{13}{16}}\left[18^{\frac{3}{4}}\max\left\{4M_{0},8L_{f}\right\}^{\frac{1}{4}}\|x_{0}-x^{*}\|\right]\epsilon^{-\frac{1}{4}}.
\label{eq:4.19}
\end{equation}
\end{theorem}

\begin{proof}
By Lemmas \ref{lem:4.2} and \ref{lem:4.3}, for $t\in\left\{0,\ldots,T\right\}$ we have
\begin{equation*}
A_{t}\tilde{f}(x_{t})\leq \varphi_{t}^{*}\leq A_{t}\tilde{f}(x^{*})+\dfrac{1}{4}\|x^{*}-x_{0}\|^{4}.
\end{equation*}
Consequently, for $t\in\left\{1,\ldots,T\right\}$,
\begin{equation}
\tilde{f}(x_{t})-\tilde{f}(x^{*})\leq\dfrac{\|x^{*}-x_{0}\|^{4}}{4A_{t}}.
\label{eq:4.20}
\end{equation}
By (\ref{eq:4.1}) and (\ref{eq:4.6}) we also have
\begin{equation*}
a_{t}^{4}=\dfrac{4^{2}(A_{t}+a_{t})^{3}}{18^{3}(2^{i_{t}}M_{t})}=\dfrac{4^{2}A_{t+1}^{3}}{18^{3}(2M_{t+1})}\geq\dfrac{4^{2}}{18^{3}\max\left\{4M_{0},8L_{f}\right\}}A_{t+1}^{3}.
\end{equation*}
Since $a_{t}=A_{t+1}-A_{t}$, it follows that
\begin{equation}
A_{t+1}-A_{t}\geq\left(\dfrac{4^{2}}{18^{3}\max\left\{4M_{0},8L_{f}\right\}}\right)^{\frac{1}{4}}A_{t+1}^{\frac{3}{4}}.
\label{eq:4.21}
\end{equation}
Let us denote
\begin{equation}
B_{t}=\dfrac{18^{3}\max\left\{4M_{0},8L_{f}\right\}}{4^{2}}A_{t}.
\label{eq:4.22}
\end{equation}
Then, by (\ref{eq:4.21}) we have
\begin{equation}
B_{t+1}-B_{t}\geq B_{t+1}^{\frac{3}{4}}.
\label{eq:4.23}
\end{equation}
Since $B_{0}=0$, it follows from (\ref{eq:4.23}) that $B_{1}\geq 1$. Thus, by Lemma A.4 in \cite{GN02} (with $\alpha=3/4$), we obtain
\begin{eqnarray*}
B_{t}&\geq & \left[\left(1-\alpha\right)\left(\dfrac{B_{1}^{\frac{1}{1-\alpha}}}{B_{1}^{\frac{1}{1-\alpha}}+1}\right)^{\alpha}\right]^{\frac{1}{1-\alpha}}(t-1)^{\frac{1}{1-\alpha}}\\
         & =    & \left[\left(\dfrac{1}{4}\right)\left(\dfrac{1}{2}\right)^{\frac{3}{4}}\right](t-1)^{4},\quad\forall t\in\left\{2,\ldots,T\right\}.
\end{eqnarray*}
Now, in view of (\ref{eq:4.22}), it follows that
\begin{equation}
A_{t}\geq\dfrac{4^{2}}{18^{3}\max\left\{4M_{0},8L_{f}\right\}2^{\frac{11}{4}}}(t-1)^{4}=\dfrac{2^{\frac{5}{4}}}{18^{3}\max\left\{4M_{0},8L_{f}\right\}}(t-1)^{4},\quad\forall t\in\left\{2,\ldots,T\right\}.
\label{eq:4.24}
\end{equation}
Then, combining (\ref{eq:4.20}) and (\ref{eq:4.24}) we get (\ref{eq:4.18}). Finally, it follows from (\ref{eq:4.17}) and (\ref{eq:4.18}) that
\begin{equation*}
\epsilon<\tilde{f}(x_{T})-\tilde{f}(x^{*})\leq\dfrac{18^{3}\max\left\{4M_{0},8L_{f}\right\}\|x_{0}-x^{*}\|^{4}}{2^{\frac{13}{4}}(T-1)^{4}},
\end{equation*}
which implies (\ref{eq:4.19}).
\end{proof}

\noindent Combining (\ref{eq:4.7}) and (\ref{eq:4.19}), we have the following worst-case complexity result. 
\begin{corollary}
\label{cor:4.1}
Under the same assumptions of Theorem \ref{thm:4.1}, let $x_{T+1}=x^{+}_{T,i_{T}}$ be the first trial point generated by Algorithm 3 such that
\begin{equation}
\min\left\{\|\nabla \tilde{f}(x_{T+1})\|,\tilde{f}(x_{T+1})-\tilde{f}(x^{*})\right\}\leq\epsilon.
\label{eq:4.25}
\end{equation}
Then, the total number of executions of Algorithm 1 is bounded as follows
\small
\begin{equation}
E_{T}\leq 4+2\left(\frac{1}{2}\right)^{\frac{13}{16}}\left[18^{\frac{3}{4}}\max\left\{4M_{0},8L_{f}\right\}^{\frac{1}{4}}\|x_{0}-x^{*}\|\right]\epsilon^{-\frac{1}{4}}+\log_{2}\left(\max\left\{2M_{0},4L_{f}\right\}\right)-\log_{2}(M_{0}).
\label{eq:4.26}
\end{equation}
\normalsize
\end{corollary}

In what follows we will assume that for all $z\in\dom\psi$,
\begin{equation}
\|g_{\psi}(z)\|\leq L_{\psi}\quad\forall g_{\psi}(z)\in\partial\psi(z).
\label{eq:4.26extra}
\end{equation}
The following lemma gives an upper bound for the number of iterations required at each execution of Algorithm 1 as inner solver in Algorithm 3.
\begin{lemma}
\label{lem:4.4}
Suppose that $f:\mathbb{R}^{n}\to\mathbb{R}$ satisfies (\ref{eq:2.2}) and that (\ref{eq:3.9}) and (\ref{eq:4.26extra}) hold. Assume that $\left\{x_{t}\right\}_{t=0}^{T}$ is generated by Algorithm 3 with
\begin{equation}
\|\nabla \tilde{f}(x^{+}_{t,i})\|>\epsilon,\quad i=0,\ldots,i_{t}\quad\text{and}\quad t=0,\ldots,T.
\label{eq:4.27}
\end{equation}
Then, for each $t\in\left\{0,\ldots,T\right\}$ and each $i\in\left\{0,\ldots,i_{t}\right\}$, the computation of $x^{+}_{t,i}$ requires at most
\begin{equation*}
\mathcal{O}\left(\log_{1.2}\left(\dfrac{3^{8}\left(7\hat{L}_{0}\right)^{4}\hat{\beta}_{0}}{4M_{0}}\epsilon^{-4}\right)\right)
\end{equation*}
iterations of Algorithm 1, where
\begin{equation}
\hat{L}_{0}=20n L_{f}R_{0}^{2}+2\text{trace}\left(\nabla^{2}f(x_{0})\right)+3\max\left\{2M_{0},4L_{f}\right\}\left(\dfrac{48}{M_{0}}\right)^{\frac{2}{3}}\hat{G}_{0}^{\frac{2}{3}},
\label{eq:4.28}
\end{equation}
\begin{equation}
\hat{\beta}_{0}=\left(10n L_{f}R_{0}^{2}+\text{trace}\left(\nabla^{2}f(x_{0})\right)\right)\left(\dfrac{48}{M_{0}}\right)^{\frac{2}{3}}\hat{G}_{0}^{\frac{2}{3}}+\dfrac{\max\left\{2M_{0},4L_{f}\right\}}{4}\left(\dfrac{48}{M_{0}}\right)^{\frac{4}{3}}\hat{G}_{0}^{\frac{4}{3}},
\label{eq:4.29}
\end{equation}
and
\begin{equation}
\hat{G}_{0}=\dfrac{4^{3}L_{f}}{6}R_{0}^{3}+4\|\nabla^{2}f(x_{0})\|R_{0}+\dfrac{4^{2}}{2}\|D^{3}f(x_{0})\|R_{0}^{2}+L_{\psi},
\label{eq:4.30}
\end{equation}
with $R_{0}$ given in (\ref{eq:3.9}) and $L_{\psi}$ given in (\ref{eq:4.26extra}).
\end{lemma}

\begin{proof}
By Theorem \ref{thm:2.1}, the computation of $x^{+}_{t,i}$ requires at most
\begin{equation*}
\mathcal{O}\left(\log_{1.2}\left(\dfrac{3^{8}\left(7L_{z_{t,i},2^{i}M_{t}}\right)^{4}\beta_{z_{t,i},2^{i}M_{t}}}{2(2^{i}M_{t})}\epsilon^{-4}\right)\right),
\end{equation*}
where 
\begin{equation}
L_{z_{t,i},2^{i}M_{t}}=\text{trace}\left(\nabla^{2}f(z_{t,i})\right)+\dfrac{3(2^{i}M_{t})}{2}\left(\dfrac{96}{2^{i}M_{t}}\right)^{\frac{2}{3}}\|\nabla\tilde{f}(z_{t,i})\|^{\frac{2}{3}}
\label{eq:4.31}
\end{equation}
and
\begin{equation}
\beta_{z_{t,i},2^{i}M_{t}}=\dfrac{1}{2}\text{trace}\left(\nabla^{2}f(z_{t,i})\right)\left(\dfrac{96}{2^{i}M_{t}}\right)^{\frac{2}{3}}\|\nabla\tilde{f}(z_{t,i})\|^{\frac{2}{3}}+\dfrac{2^{i}M_{t}}{8}\left(\dfrac{96}{2^{i}M_{t}}\right)^{\frac{4}{3}}\|\nabla\tilde{f}(z_{t,i})\|^{\frac{4}{3}}.
\label{eq:4.32}
\end{equation}
It follows from Step 1 of Algorithm 3 and Lemma \ref{lem:4.1} that
\begin{equation}
2M_{0}\leq 2^{i}M_{t}\leq 2^{i_{t}}M_{t}=2M_{t+1}\leq 2\max\left\{2M_{0},4L_{f}\right\}.
\label{eq:4.33}
\end{equation}
Thus, the computation of $x^{+}_{t,i}$ requires at most
\begin{equation*}
\mathcal{O}\left(\log_{1.2}\left(\dfrac{3^{8}\left(7L_{z_{t,i},2^{i}M_{t}}\right)^{4}\beta_{z_{t,i},2^{i}M_{t}}}{4M_{0}}\epsilon^{-4}\right)\right)
\end{equation*}
Now, it is enough to show that $L_{z_{t,i},2^{i}M_{t}}$ and $\beta_{z_{t,i},2^{i}M_{t}}$ are bounded from above by $\hat{L}_{0}$ and $\hat{\beta}_{0}$, respectively. In view of (\ref{eq:2.6}) we have
\begin{eqnarray*}
\nabla^{2}f(z_{t,i})&=&\nabla^{2}f(z_{t,i})-\nabla\Phi_{x_{0}}(z_{t,i})+\nabla^{2}\Phi_{x_{0}}(z_{t,i})\\
                               &\preceq & \dfrac{L_{f}}{2}\|z_{t,i}-x_{0}\|^{2}I+\nabla^{2}f(x_{0})+D^{3}f(x_{0})[z_{t,i}-x_{0}],
\end{eqnarray*}
and so
\begin{equation}
\text{trace}\left(\nabla^{2}f(z_{t,i})\right)\leq\dfrac{nL_{f}}{2}\|z_{t,i}-x_{0}\|^{2}+\text{trace}\left(\nabla^{2}f(x_{0})\right)+\text{trace}\left(D^{3}f(x_{0})[z_{t,i}-x_{0}]\right).
\label{eq:4.34}
\end{equation}
On the other hand, by Lemma 3 in \cite{NES1} we also have
\begin{equation*}
D^{3}f(x_{0})[z_{t,i}-x_{0}]\preceq\nabla^{2}f(x_{0})+\dfrac{L_{f}}{2}\|z_{t,i}-x_{0}\|^{2}I,
\end{equation*}
which gives
\begin{equation}
\text{trace}\left(D^{3}f(x_{0})[z_{t,i}-x_{0}]\right)\leq\text{trace}\left(\nabla^{2}f(x_{0})\right)+\dfrac{n L_{f}}{2}\|z_{t,i}-x_{0}\|^{2}.
\label{eq:4.35}
\end{equation}
Thus, combining (\ref{eq:4.34}) and (\ref{eq:4.35}) we obtain
\begin{eqnarray}
\text{trace}\left(\nabla^{2}f(z_{t,i})\right)&\leq & n L_{f}\|z_{t,i}-x_{0}\|^{2}+2\text{trace}\left(\nabla^{2}f(x_{0})\right)\nonumber\\
                                                                    &\leq & n L_{f}\left(\|z_{t,i}-x^{*}\|+\|x^{*}-x_{0}\|\right)^{2}+2\text{trace}\left(\nabla^{2}f(x_{0})\right)\nonumber\\
                                                                    &\leq & 2n L_{f}\left(\|z_{t,i}-x^{*}\|^{2}+\|x_{0}-x^{*}\|^{2}\right)+2\text{trace}\left(\nabla^{2}f(x_{0})\right).
\label{eq:4.36}
\end{eqnarray}
By (\ref{eq:4.2}) we have
\begin{eqnarray}
\|z_{t,i}-x^{*}\|&=&\|(1-\gamma_{t,i})x_{t}+\gamma_{t,i}v_{t}-(1-\gamma_{t,i})x^{*}-\gamma_{t,i}x^{*}\|\nonumber\\
                           &\leq & (1-\gamma_{t,i})\|x_{t}-x^{*}\|+\gamma_{t,i}\|v_{t}-x^{*}\|\nonumber\\
                           &\leq & \|x_{t}-x^{*}\|+\|v_{t}-x^{*}\|,
\label{eq:4.37}
\end{eqnarray}
where the last inequality is due to $\gamma_{t,i}\in (0,1]$. Since $\varphi_{t}(\,\cdot\,)$ is uniformly convex of degree $4$ with parameter $1/4$, by Lemmas \ref{lem:4.2} and \ref{lem:4.3} we also have
\begin{equation*}
\dfrac{1}{16}\|v_{t}-x\|^{4}+A_{t}\tilde{f}(x_{t})\leq\dfrac{1}{16}\|v_{t}-x\|^{4}+\varphi_{t}^{*}\leq\varphi_{t}(x)\leq A_{t}\tilde{f}(x)+\dfrac{1}{4}\|x-x_{0}\|^{4}.
\end{equation*}
In particular, for $x=x^{*}$ it follows that
\begin{equation*}
\dfrac{1}{16}\|v_{t}-x^{*}\|^{4}+A_{t}(\tilde{f}(x_{0})-\tilde{f}(x^{*}))\leq\dfrac{1}{4}\|x^{*}-x_{0}\|^{4},
\end{equation*}
and so
\begin{equation}
\|v_{t}-x^{*}\|\leq 4^{\frac{1}{4}}\|x_{0}-x^{*}\|\leq 2\|x_{0}-x^{*}\|.
\label{eq:4.38}
\end{equation}
Moreover, by Lemmas \ref{lem:4.2} and \ref{lem:4.3} we also have
\begin{equation*}
A_{t}\tilde{f}(x_{t})\leq\varphi_{t}^{*}\leq\varphi_{t}(x_{0})\leq A_{t}\tilde{f}(x_{0}).
\end{equation*}
Consequently, $x_{t}\in\mathcal{L}(x_{0})$, and by (\ref{eq:3.9}) we get
\begin{equation}
\|x_{t}-x^{*}\|\leq R_{0}.
\label{eq:4.39}
\end{equation}
Thus, combining (\ref{eq:4.37})-(\ref{eq:4.39}) and (\ref{eq:3.9}) it follows that
\begin{equation}
\|z_{t,i}-x^{*}\|\leq 3R_{0}.
\label{eq:4.40}
\end{equation}
Then, by (\ref{eq:4.36}), (\ref{eq:4.40}) and (\ref{eq:3.9}) we obtain
\begin{equation}
\text{trace}\left(\nabla^{2}f(z_{t,i})\right)\leq 20n L_{f}R_{0}^{2}+2\text{trace}\left(\nabla^{2}f(x_{0})\right)
\label{eq:4.41}
\end{equation} 
Moreover, by (\ref{eq:2.5}), (\ref{eq:4.40}) and (\ref{eq:3.9}) we have
\begin{eqnarray*}
\|\nabla f(z_{t,i})\|&=&\|\nabla f(z_{t,i})-\nabla\Phi_{x_{0}}(z_{t,i})\|+\|\nabla\Phi_{x_{0}}(z_{t,i})\|\\
                              &\leq &\dfrac{L_{f}}{6}\|z_{t,i}-x_{0}\|^{3}+\|\nabla f(x_{0})\|\|z_{t,i}-x_{0}\|+\dfrac{1}{2}\|D^{3}f(x_{0})\|\|z_{t,i}-x_{0}\|^{2}\\
                             &\leq &\dfrac{4^{3}L_{f}}{6}R_{0}^{3}+4\|\nabla f(x_{0})\|R_{0}+\dfrac{4^{2}}{2}\|D^{3}f(x_{0})\|R_{0}^{2}.
\end{eqnarray*}
Thus, by (\ref{eq:4.26extra}) and (\ref{eq:4.30}) we get
\begin{equation}
\|\nabla\tilde{f}(z_{t,i})\|\leq\hat{G}_{0}.
\label{eq:4.42}
\end{equation}
Finally, it follows from (\ref{eq:4.31}), (\ref{eq:4.41}), (\ref{eq:4.42}) and (\ref{eq:4.28}) that $L_{z_{t,i},2^{i}M_{t}}\leq \hat{L}_{0}$, while (\ref{eq:4.32}), (\ref{eq:4.41}), (\ref{eq:4.42}) and (\ref{eq:4.29}) imply that $\beta_{z_{t,i},2^{i}M_{t}}\leq \hat{\beta}_{0}$.
\end{proof}

\begin{remark}
In view of Corollary \ref{cor:4.1} and Lemma \ref{lem:4.4}, to generate an $\epsilon$-approximate minimizer of $\tilde{f}(\,\cdot\,)$, Algorithm 3 needs at most $\mathcal{O}\left(|\log_{1.2}(\epsilon)|\epsilon^{-\frac{1}{4}}\right)$ iterations of Algorithm 1.
\end{remark}

\section{Illustrative Numerical Results}

To investigate the effectiveness of our new adaptive strategy, we tested a MATLAB implementation of Algorithm 2. We applied the corresponding code to smooth logistic regression problems of the form
\begin{equation*}
\min_{x\in\mathbb{R}^{n+1}}\,f(x)=-\sum_{i=1}^{m}\left[b^{(i)}\log(m_{x}(a^{(i)}))+(1-b^{(i)})\log(1-m_{x}(a^{(i)}))\right],
\end{equation*}
where $\left\{(a^{(i)},b^{(i)})\right\}_{i=1}^{m}\subset\mathbb{R}^{n+1}\times\left\{0,1\right\}$ is the dataset (with $a_{1}^{(i)}=1$ for $i=1,\ldots,m$) $m_{x}(a)\equiv 1/\left(1+e^{-\langle a,x\rangle}\right)$, and $\log(\,\cdot\,)$ denotes the natural logarithm. Three datasets from \cite{DUA} were considered. For each of them we used the starting point $x_{0}=[1\,1\,\ldots \,1]^{T}$, the parameter $M_{0}=1$ and the stopping criterion $\|\nabla f(x_{k})\|\leq\epsilon$. The results are presented in Tables \ref{table:1}, \ref{table:2} and \ref{table:3}, where ``IT'' represents the number of outer iterations (i.e., iterations of Algorithm 2), ``CO'' represents the number of calls of the oracle\footnote{Each call of the oracle corresponds to either one function evaluation or one (high-order) derivative evaluation.}, ``BGM-E'' represents the number of executions of Algorithm 1, ``BGM-IT'' represents the total number of inner iterations (i.e., iterations of Algorithm 1), and ``BGM-A'' represents the average number of inner iterations per execution of Algorithm 1. 

\begin{table}[h!]
\centering
\small{    
\begin{tabular}{ | c | c | c | c | c | c | }
    \hline
    \textbf{$\epsilon$}   & IT & CO& BGM-E&  BGM-IT & BGM-A\\ \hline
		$10^{-2}$         & 42      & 252	  & 83    & 469 & 5.6506 \\\hline
		$10^{-4}$         & 42      & 252      & 83 	  & 491 &  5.9156\\\hline
		$10^{-6}$         & 43      & 256      & 84    & 496 &  5.9048\\ \hline
                      $10^{-8}$        &  43      & 256     & 85    & 520 &  6.1905\\
\hline
    \end{tabular}
}
\caption{Results for the dataset Diabetes \cite{FIS} ($n=8$ and $m=768$)}
\label{table:1}
\end{table}

\begin{table}[h!]
\centering
\small{    
\begin{tabular}{ | c | c | c | c | c | c | }
    \hline
    \textbf{$\epsilon$}   & IT & CO& BGM-E&  BGM-IT & BGM-A\\ \hline
		$10^{-2}$         & 4      & 20     & 5    & 62 &  12.4000 \\\hline
		$10^{-4}$         & 4      & 20     & 5    & 84 &  16.8000\\\hline
		$10^{-6}$         & 5      & 20     & 5    & 108 &  21.6000\\ \hline
                      $10^{-8}$        & 5      & 24     & 6    & 102 &  17.0000\\
\hline
    \end{tabular}
}
\caption{Results for the dataset Phishing \cite{ABDE} ($n=9$ and $m=1353$)}
\label{table:2}
\end{table}

\begin{table}[h!]
\centering
\small{    
\begin{tabular}{ | c | c | c | c | c | c | }
    \hline
    \textbf{$\epsilon$}   & IT & CO& BGM-E&  BGM-IT & BGM-A\\ \hline
		$10^{-2}$  &  59 &      239      & 60	   & 258     & 4.3000  \\\hline
		$10^{-4}$  &125&       503      & 126        & 522      & 4.1428 \\\hline
		$10^{-6}$  &411&       1647     & 412       & 1666    & 4.0437\\ \hline
                     $10^{-8}$   &1731&     6927      & 1732    & 6946    & 4.0104\\
\hline
    \end{tabular}
}
\caption{Results for the dataset Ionosphere \cite{SIG} ($n=34$ and $m=351$)}
\label{table:3}
\end{table}

As we can see, the new adaptive scheme is able to find approximate stationary points with efficiency as predicted by the theory. Moreover, the average number of inner iterations required to compute a suitable trial point is usually small.

\section{Second-order Implementations}

In our adaptive third-order methods, third-order information about $f(\,\cdot\,)$ is only required for the computation of the gradients
\begin{equation}
\nabla\Omega_{x,M}(y)=\nabla f(x)+\nabla^{2}f(x)(y-x)+\dfrac{1}{3}D^{3}f(x)[y-x]^{2}+\dfrac{M}{2}\|y-x\|^{2}(y-x),
\label{eq:stela6.1}
\end{equation}
which are used in the inner solver (Algorithm 2). Notice that, inequality (\ref{eq:2.5}) can be rewritten as 
\begin{equation*}
\|\nabla f(y)-\nabla f(x)-\nabla^{2}f(x)(y-x)-\dfrac{1}{2}D^{3}f(x)[y-x]^{2}\|\leq\dfrac{L_{f}}{3!}\|y-x\|^{3},\quad\forall x,y\in\mathbb{R}^{n}.
\end{equation*}
In particular, for any $h\in\mathbb{R}^{n}$ and $\tau>0$, we have
\begin{equation*}
\left\|\nabla f(x+\tau h)-\nabla f(x)-\tau\nabla^{2}f(x)h-\frac{\tau^{2}}{2}D^{3}f(x)[h]^{2}\right\|\leq\dfrac{L_{f}}{3!}\|h\|^{3}
\end{equation*}
and
\begin{equation*}
\left\|\nabla f(x-\tau h)-\nabla f(x)+\tau\nabla^{2}f(x)h-\frac{\tau^{2}}{2}D^{3}f(x)[h]^{2}\right\|\leq\dfrac{L_{f}}{3!}\|h\|^{3}
\end{equation*}
Combining these two inequalities, it follows that
\begin{equation*}
\|\nabla f(x+\tau h)+\nabla f(x-\tau h)-2\nabla f(x)-\tau^{2}D^{3}f(x)[h]^{2}\|\leq\dfrac{L_{f}}{3}\tau^{3}\|h\|^{3},
\end{equation*}
and so we can approximate $D^{3}f(x)[h]^{2}$ by the vector
\begin{equation}
T_{\tau}(h)=\dfrac{\nabla f(x+\tau h)+\nabla f(x-\tau h)-2\nabla f(x)}{\tau^{2}}
\label{eq:stela6.2}
\end{equation}
with an error bounded as follows:
\begin{equation*}
\|T_{\tau}(h)-D^{3}f(x)[h]^{2}\|\leq\dfrac{L_{f}}{3}\tau\|h\|^{3}.
\end{equation*}
Thus, by replacing $D^{3}f(x)[y-x]^{2}$ in (\ref{eq:stela6.1}) by $T_{\tau}(y-x)$, we obtain the vector
\begin{equation}
g_{x,M,\tau}(y)=\nabla f(x)+\nabla^{2}f(x)(y-x)+\dfrac{1}{3}T_{\tau}(y-x)+\dfrac{M}{2}\|y-x\|^{2}(y-x),
\label{eq:stela4.3}
\end{equation}
which approximates $\nabla\Omega_{x,M}(y)$ with an error bound
\begin{equation*}
\|g_{x,M,\tau}(y)-\nabla\Omega_{x,M}(y)\|\leq\dfrac{L_{f}}{9}\tau\|y-x\|^{3}.
\end{equation*}
Since $g_{x,M,\tau}(\,\cdot\,)$ defined in (\ref{eq:stela4.3}) depends only on first and second-order information about $f(\,\cdot\,)$, this remark opens the possibility for the development of second-order implementations of our third-order methods. As shown in \cite{NES5}, with a carefull choice of $\tau$, it is indeed possible to obtain this type of second-order scheme preserving the global convergence rate of the original third-order method. However, the choice for $\tau$ in \cite{NES5} requires the knowledge of the Lipschitz constant $L_{f}$. Therefore, the development of adaptive strategies in this context remains a very interesting open question.
\section{Conclusion}

In this paper we presented adaptive third-order methods for composite convex optimization problems in which the smooth part is a three-times continuously differentiable function with Lipschitz continuous third-order derivatives. Trial points are computed using a Bregman gradient method as inner solver. Our adaptive procedure to update the regularization parameters in the third-order tensor models take into account the complexity of minimizing these models. When the regularization parameter is sufficiently large, the inner solver is guaranteed to find a suitable approximate stationary point of the model within $\mathcal{O}\left(|\log(\epsilon)|\right)$ iterations. Thus, once a slow convergence rate is detected, the execution of the inner solver is stopped and the regularization parameter is increased.  Using this strategy we obtained a basic adaptive third-order method that finds an $\epsilon$-approximate minimizer of the objective function performing at most $\mathcal{O}\left(|\log(\epsilon)|\epsilon^{-\frac{1}{3}}\right)$ iterations of the inner solver. We also obtained an accelerated adaptive third-order method with an improved complexity bound of $\mathcal{O}\left(|\log(\epsilon)|\epsilon^{-\frac{1}{4}}\right)$ inner iterations. Preliminary numerical experiments illustrated the effectiveness of our new adaptive technique.


\end{document}